\documentclass[12pt,journal, onecolumn, draftcls]{IEEEtran}

\usepackage{verbatim}%comment
\usepackage{amsmath,graphicx}
\usepackage{mathrsfs}
\usepackage{cite}%bibliography
\usepackage{amsthm}%theorenm, lemma, definition, ...
\usepackage{epsfig,psfrag}%subfigure,enumerate
\usepackage{amsfonts,amssymb}%math fonts and symbols

\newcommand{\bR}{\boldsymbol{R}}
\newcommand{\bP}{\boldsymbol{P}}
\newcommand{\ba}{\boldsymbol{a}}
\newcommand{\bx}{\boldsymbol{x}}
\newcommand{\bn}{\boldsymbol{n}}
\newcommand{\bs}{\boldsymbol{s}}
\newcommand{\bA}{\boldsymbol{A}}
\newcommand{\bI}{\boldsymbol{I}}
\newcommand{\bS}{\boldsymbol{S}}
\newcommand{\bg}{\boldsymbol{g}}
\newcommand{\be}{\boldsymbol{e}}
\newcommand{\bG}{\boldsymbol{G}}
\newcommand{\bE}{\boldsymbol{E}}

\newcommand{\bQ}{\boldsymbol{Q}}

\newcommand{\bT}{\boldsymbol{T}}
\newcommand{\bV}{\boldsymbol{V}}
\newcommand{\bB}{\boldsymbol{B}}

\begin{document}
\title{Subspace Leakage Analysis and Improved DOA Estimation with
Small Sample Size}

\author{Mahdi~Shaghaghi and~Sergiy~A.~Vorobyov
\thanks{M. Shaghaghi is with the Department
of Electrical and Computer Engineering, University of Alberta,
Edmonton, AB, T6G 2V4 Canada (e-mail: mahdi.shaghaghi@ualberta.ca).
S.~A.~Vorobyov is with Aalto University, Department of Signal
Processing and
Acoustics, Finland (e-mail: sergiy.vorobyov@aalto.fi). S.~A.~Vorobyov
is the corresponding author.}
\thanks{Parts of this paper have been/will be published at the IEEE
Inter. Workshop on Computational Advances in Multi-Sensor Adaptive
Processing (CAMSAP), The Friendly Island, Saint Martin, Dec.~2013
and the IEEE Inter. Conf. Acoustics, Speech, and Signal Processing
(ICASSP),  Brisbane, Australia, Apr.~2015.}
\thanks{This work was supported in part by the Natural Sciences and
Engineering Research Council (NSERC) of Canada.}
}

%\markboth{IEEE TRANSACTIONS ON SIGNAL PROCESSING}%
%{Shell \MakeLowercase{\textit{et al.}}: Bare Demo of IEEEtran.cls
%for Journals}
\maketitle

\begin{abstract}
Classical methods of DOA estimation such as the MUSIC algorithm are
based on estimating the signal and noise subspaces from the sample
covariance matrix. For a small number of samples, such methods are
exposed to performance breakdown, as the sample covariance matrix
can largely deviate from the true covariance matrix. In this paper,
the problem of DOA estimation performance breakdown is investigated.
We consider the structure of the sample covariance matrix and the
dynamics of the root-MUSIC algorithm. The performance breakdown in
the threshold region is associated with the subspace leakage where
some portion of the true signal subspace resides in the estimated
noise subspace. In this paper, the subspace leakage is theoretically
derived. We also propose a two-step method which improves the
performance by modifying the sample covariance matrix such that the
amount of the subspace leakage is reduced. Furthermore, we introduce
a phenomenon named as root-swap which occurs in the root-MUSIC
algorithm in the low sample size region and degrades the performance
of the DOA estimation. A new method is then proposed to alleviate
this problem. Numerical examples and simulation results are given
for uncorrelated and correlated sources to illustrate the
improvement achieved by the proposed methods. Moreover, the proposed
algorithms are combined with the pseudo-noise resampling method to
further improve the performance.
\end{abstract}


\IEEEpeerreviewmaketitle

\section{Introduction}
%\IEEEPARstart{T}{his}
Classical parameter estimation methods of direction-of-arrival
(DOA), Doppler shifts, frequencies, time delays, etc. such as the
multiple signal classification (MUSIC) \cite{Schmidt86:MUSIC},
root-MUSIC \cite{Barabell83:root_Music}, and estimation of signal
parameters via rotational invariance techniques (ESPRIT)
\cite{Roy89:ESPRIT} are based on estimating the signal and noise
subspaces from the sample data covariance matrix. It is well-known
that these methods suffer from performance breakdown for a small
number of samples or low signal-to-noise ratio (SNR) values where
the expected estimation error departs from the Cram\'{e}r-Rao bound
(CRB)\cite{Scharf95:prob_subspace_swap}. The SNR region at which
this phenomenon happens is known as the threshold region.

The fidelity of the sample data covariance matrix to the true data
covariance matrix plays a critical role in a successful estimation.
At the low SNR and/or small sample size region, the sample data
covariance matrix can largely deviate from the true one. There are
various methods introduced in the literature which target at
improving the estimation of the covariance matrix
\cite{Carlson88:diagonal_loading,
Eldar10:Shrinkage_MMSE,Pillai89:Forward_Backward,
Evans82:Spatial_Smoothing,Mestre08:TSP_DOA,
Gershman97:Pseudo_random_DOA, Vasylyshyn13:Pseudo_noise_rootMusic,
Qian14:impro_uni_r_Music_Pseudo_noise}.

Diagonal loading \cite{Carlson88:diagonal_loading} and
shrinkage-based \cite{Eldar10:Shrinkage_MMSE} methods improve the
estimate of the data covariance matrix by scaling and shifting the
eigenvalues of the sample data covariance matrix. However, the
eigenvectors are kept unchanged. As a result, the estimated signal
and noise projection matrices from the improved covariance matrices
are exactly the same as those obtained from the sample data
covariance matrix. Therefore, these methods are not really
beneficial for the subspace-based parameter estimation algorithms.

Data covariance matrix estimation can be also improved by the
means of using forward-backward averaging (FB)
\cite{Pillai89:Forward_Backward} and spatial smoothing-based
techniques \cite{Evans82:Spatial_Smoothing}. The effect of FB is
known to be equivalent to approximately doubling the number of
samples. Thus, the covariance estimate improves accordingly. The
spatial smoothing technique can also be interpreted as virtually
increasing the number of samples at the cost of averaging over
sub-arrays of smaller size compared to the whole array.
These approaches can also decorrelate pairs (in case of FB) or more
correlated source signals.
In \cite{Mestre08:TSP_DOA}, techniques from random matrix theory
have been developed to improve the performance of the MUSIC
algorithm. The introduced method considers the asymptotic situation
when both the sample size and the number of array elements tend to
infinity at the same rate. It is then inferred that the improved
method gives a more accurate description of the situation when these
two quantities are finite and comparable in magnitude. However, the
performance of the introduced method is not satisfactory at the
small sample size scenario \cite{MSH_13:i_root_Music}.

A more promising approach to remedy the performance breakdown at the
threshold region was introduced in
\cite{Gershman97:Pseudo_random_DOA} and has been further improved in
\cite{Vasylyshyn13:Pseudo_noise_rootMusic} and
\cite{Qian14:impro_uni_r_Music_Pseudo_noise}. These methods are
based on a technique called pseudo-noise resampling which uses
synthetically generated pseudo-noise to perturb the original noise.
The pseudo-noise is added to the observed data, and a new estimate
of the covariance matrix is obtained, which leads to new DOA
estimates. This process is repeated for a number of times, and the
final DOAs are determined based on the bank of the DOA estimates.

In this paper, we tackle the problem of the performance breakdown at
the threshold region by considering the structure of the sample data
covariance matrix and the dynamics of the root-MUSIC algorithm. It
is shown in \cite{Mestre08:TSP_sub_leak} that the performance
breakdown problem is associated with the inter-subspace leakage
``whereby a small portion of the true signal eigenvector resides in
the sample noise subspace (and vice versa)''. The notion of leakage
comes originally from the performance assessment strategy based on
the first order approximation of the estimation error caused by the
perturbed subspace estimate, which happens because of the additive
noise contribution
\cite{Li93:Perform_unification,Xu02:Perturbation_Subspace,
Liu08:Perturbation_SVD,Haart14:ESPRIT_leakage}. This approach
directly models the leakage of the noise subspace into signal
subspace and allows to compute the corresponding perturbation matrix
between the components of the subspaces. Here, we formally define
the \textit{subspace leakage} notion as a Frobenius norm of the
perturbation matrix, and we present its theoretical derivation. We
propose a two-step method which improves the performance of the
root-MUSIC algorithm by modifying the sample data covariance matrix
such that the amount of the subspace leakage is reduced.
Furthermore, we introduce a phenomenon named as \textit{root-swap}
which occurs in the root-MUSIC algorithm at the threshold region and
degrades the performance of the parameter estimation. A new method
is then proposed to alleviate this problem.

It will be shown that there are undesirable by-products in the
sample data covariance matrix that tend to zero as the number of
samples goes to infinity. However, for a limited number of samples,
these terms can have significant values leading to a large amount of
subspace leakage. One possible approach to remedy the effect of the
undesirable components is to consider the \textit{eigenvalue
perturbation} caused by these terms. The incorporation of this
knowledge into the estimation method can result in better estimates
of the signal and noise subspaces. In this paper, we propose a
two-step algorithm in order to reduce the effect of the undesirable
terms. The introduced method is based on estimating the parameters
at the first step and modifying the covariance matrix using the
estimated parameters at the second step. We will theoretically
derive the subspace leakage at both steps. Then, it will be shown
using numerical examples that the subspace leakage is reduced at the
second step leading to better performance.

In the root-MUSIC method, the estimation error of the roots has a
variance which is proportional to the variance of noise over the
number of samples \cite{Proakis92:root_Music_performance}.
Therefore, at the threshold region, the variance of the estimation
error can have a significant value which in turn can result in a
swap between a root corresponding to a signal source with another
root which is not associated with any signal source. We dub this
phenomenon as root-swap. Then, a new method is proposed to remedy
this problem. The introduced method considers different combinations
of the roots as the candidates for the signal sources. These
candidates are then evaluated using the stochastic maximum
likelihood (SML) function, and the combination that minimizes the
objective function is picked up for the parameter estimates.

The rest of the paper is organized as follows. The system model is
given and the root-MUSIC algorithm is briefly reviewed in
Section~\ref{sec:Sys_Mod_Sub_Leak}. The two-step and root-swap
algorithms are proposed in
Section~\ref{sec:Proposed_2_step_root_swap}. Subspace leakage is
defined and theoretically derived in Section~\ref{sec:Sub_leak}.
Numerical examples and simulation results are given in
Section~\ref{sec:Sim_Sub_Leak}. Section~\ref{sec:conclude} concludes
the paper. Appendix~\ref{sec:appndx_root_swap_prob} gives an
approximation for the probability of root-swap, and finally, the
details of the subspace leakage derivation for the two-step
root-MUSIC algorithm are presented in
Appendices~\ref{sec:appndx_SL1} and \ref{sec:appndx_SL2}.

\section{System Model and Background}
\label{sec:Sys_Mod_Sub_Leak}
An example of a noise-corrupted linear superposition of $K$ undamped
exponentials received by $M$ ($M > K$) antennas is the array
processing model. Thus, consider $K$ number of narrowband plane
waves impinging on a uniform linear array (ULA) from directions
$\theta_1,~\theta_2,\cdots,~\theta_K$. Without loss of generality,
assume $-\pi/2 \leq \theta_1 \leq \theta_2 \leq \cdots \leq \theta_K
\leq \pi/2$. The antenna elements are separated from each other by a
distance of $d \leq \lambda / 2$ where $\lambda$ is the wavelength
of the plane wave impinging on the array. The steering vector of the
array $\ba(\theta) \in \mathbb{C}^{M \times 1}$ is then given as
\begin{equation}
\ba(\theta) \triangleq \left[1,~e^{-j2 \pi (d/\lambda)
\sin(\theta)},\cdots,~e^{-j2 \pi (M-1)(d/\lambda)
\sin(\theta)}\right]^T \label{eq:steering_vec}
\end{equation}
where $\left(\cdot\right)^T$ stands for the transposition operator.
At time instant $t \in \mathbb{N}$, the received vector $\bx(t) \in
\mathbb{C}^{M \times 1}$ is given by
\begin{equation}
\bx(t) = \sum_{i = 1}^K \ba(\theta_i) s_i(t) + \bn(t)
\label{eq:sys_mod_1}
\end{equation}
where $s_i(t) \in \mathbb{C}$ is the amplitude of the $i$-th wave
(source) and $\bn(t) \in \mathbb{C}^{M \times 1}$ is the noise
vector at time $t$. By arranging the amplitudes of the sources in
the vector $\bs(t) = \left[s_1(t),s_2(t),\cdots,s_K(t)\right]^T \in
\mathbb{C}^{K \times 1}$ and forming the Vandermonde matrix $\bA =
\left[\ba(\theta_1),~\ba(\theta_2),\cdots,~\ba(\theta_K)\right] \in
\mathbb{C}^{M \times K}$, the model \eqref{eq:sys_mod_1} can be
rewritten in matrix-vector form as
\begin{equation}
\bx(t) = \bA \bs(t) + \bn(t). \label{eq:sys_mod_2}
\end{equation}
We consider the noise vector $\bn(t)$ to be independent from the
sources and noise vectors at other time instances and to have the
circularly-symmetric complex jointly-Gaussian distribution
$\mathcal{N}_C(0,\sigma_{\text{n}}^2\bI_M)$ where $\bI_M$ is the
identity matrix of size $M$. Considering the system model
\eqref{eq:sys_mod_2}, the data covariance matrix $\bR \in
\mathbb{C}^{M \times M}$ is given by
\begin{equation}
\bR \triangleq E\left\{\bx(t) \bx^H(t) \right\} = \bA\bS\bA^H +
\sigma_{\text{n}}^2\bI_M \label{eq:R_model}
\end{equation}
where $\bS = E \left\{\bs(t) \bs^H(t) \right\} \in \mathbb{C}^{K
\times K}$ is the source covariance matrix and
$\left(\cdot\right)^H$ and $E\{\cdot\}$ stand for the Hermitian
transposition and the expectation operators, respectively.

Let $\lambda_1 \leq \lambda_2 \leq \cdots \leq \lambda_M$ be the
eigenvalues of $\bR$ arranged in nondecreasing order, and let
$\bg_1,~\bg_2,\cdots,~\bg_{M-K}$ be the noise eigenvectors
associated with $\lambda_1,~\lambda_2,\cdots,~\lambda_{M-K}$ and
$\be_1,~\be_2,\cdots,~\be_K$ be the signal eigenvectors
corresponding to
$\lambda_{M-K+1},~\lambda_{M-K+2},\cdots,~\lambda_M$. Let also $\bG
\in \mathbb{C}^{M \times (M-K)}$ and $\bE \in \mathbb{C}^{M \times
K}$ be defined as $\bG \triangleq
\left[\bg_1,~\bg_2,\cdots,~\bg_{M-K}\right]$ and $\bE \triangleq
\left[\be_1,~\be_2,\cdots,~\be_K\right]$. The range spaces of $\bG$
and $\bE$ represent the true noise and signal subspaces,
respectively.

Let $N$ number of snapshots (samples) be available. The basic method
for estimating the data covariance matrix from the samples $\bx(t)$
($1 \leq t\leq N$) is
\begin{equation}
\widehat{\bR} \triangleq \frac{1}{N}\sum_{t=1}^N \bx(t) \bx^H(t)
\label{eq:Rhat_conventional}
\end{equation}
where $\widehat{\bR} \in \mathbb{C}^{M \times M}$ is the sample data
covariance matrix. Consider the eigendecomposition of
$\widehat{\bR}$. Let
$\hat{\bg}_1,~\hat{\bg}_2,\cdots,~\hat{\bg}_{M-K}$ be the estimated
noise eigenvectors and
$\hat{\be}_1,~\hat{\be}_2,\cdots,~\hat{\be}_K$ be the estimated
signal eigenvectors. Form $\widehat{\bG} \in \mathbb{C}^{M \times
(M-K)}$ and $\widehat{\bE} \in \mathbb{C}^{M \times K}$ by placing
the estimated noise and signal eigenvectors as the columns of
$\widehat{\bG}$ and $\widehat{\bE}$, respectively. The range spaces
of $\widehat{\bG}$ and $\widehat{\bE}$ represent the estimations of
the noise and signal subspaces, respectively.

Recalling \eqref{eq:steering_vec} and defining $z \triangleq e^{j2
\pi (d/\lambda) \sin(\theta)}$, the steering vector can be rewritten
as $\ba(z) =
\left[1,z^{-1},\cdots,z^{-(M-1)}\right]^T\,\hspace{-2mm}$. In the
root-MUSIC method, the roots of the equation $\ba^T(z^{-1})
\widehat{\bG} \widehat{\bG}^H \hspace{-1mm} \ba(z)\hspace{-1mm} = 0$
which are located inside the unit circle are considered. These roots
are sorted based on their distance to the unit circle, and the first
$K$ number of the roots which are closer to the unit circle are
picked. The estimates of the DOAs denoted by
$\hat{\theta}_1,~\hat{\theta}_2,\cdots,~\hat{\theta}_K$ are then
obtained by multiplying the angles of the selected roots by
$\lambda/(2 \pi d)$ and taking the inverse sinusoid function of the
results.

%Next, the eigenvalues of $\widehat{\bR}$ are computed and
%sorted in nondecreasing order. The eigenvectors corresponding to the
%first $(M-K)$ eigenvalues are then arranged as the columns of a
%matrix denoted by $\widehat{\bG} \in \mathbb{C}^{M \times
%(M-K)}$. The range space of $\widehat{\bG}$ represents
%the estimation of the noise subspace. Recalling
%\eqref{eq:steering_vec} and defining $z \triangleq e^{j2 \pi
%(d/\lambda) \sin(\theta)}$, the steering vector can be rewritten as
%\begin{equation}
%\ba(z) = \left[1,~z^{-1},\cdots,~z^{-(M-1)}\right]^T.
%\end{equation}
%For the root-MUSIC method, the roots of the equation
%$\ba^T(z^{-1}) \widehat{\bG}
%\widehat{\bG}^H \ba(z) = 0$ are used to obtain
%the estimates of the DOAs denoted by
%$\hat{\theta}_1,~\hat{\theta}_2,\cdots,~\hat{\theta}_K$
%\cite{Barabell83:root_Music}.

\section{Proposed Methods}
\label{sec:Proposed_2_step_root_swap}

\subsection{Two-step root-MUSIC algorithm}

Let us start by expanding \eqref{eq:Rhat_conventional} using
\eqref{eq:sys_mod_2} as follows
\begin{eqnarray}
\hspace{-3mm}\widehat{\bR} \hspace{-3mm} &=& \hspace{-3mm}
\frac{1}{N}\sum_{t=1}^N \left( \bA \bs(t) + \bn(t) \right)
\left(\bA \bs(t) + \bn(t)\right)^H \nonumber \\
\hspace{-3mm} \hspace{-3mm} &=& \hspace{-3mm} \bA \left\{
\frac{1}{N}\sum_{t=1}^N \bs(t) \bs^H(t) \right\} \bA^H +
\frac{1}{N}\sum_{t=1}^N
\bn(t) \bn^H(t) \nonumber \\
\hspace{-3mm} \hspace{-3mm} && \hspace{-3mm} + \bA \left\{
\frac{1}{N}\sum_{t=1}^N \bs(t) \bn^H(t) \right\} + \left\{
\frac{1}{N}\sum_{t=1}^N \bn(t) \bs^H(t) \right\} \bA^H.
\label{eq:Rhat_expand_2}
\end{eqnarray}
Comparing \eqref{eq:Rhat_expand_2} with \eqref{eq:R_model}, it can
be observed that the expansion of $\widehat{\bR}$ consists of four
terms while the model for $\bR$ comprises two summands. The first
two terms of $\widehat{\bR}$ given by \eqref{eq:Rhat_expand_2} can
be considered as estimates for the two summands of $\bR$, which
represent the signal and noise components, respectively. The last
two terms of $\widehat{\bR}$ in \eqref{eq:Rhat_expand_2} are
undesirable by-products which can be viewed as estimates for the
correlation between the signal and noise vectors. In the system
model under study, we consider the noise vectors to be zero-mean and
also independent of the signal vectors. Therefore, the signal and
noise components are uncorrelated to each other. As a result, for a
large enough number of samples $N$, the last two terms in
\eqref{eq:Rhat_expand_2} tend to zero. However, the number of
available samples can be limited in practical applications. In this
case, the last two terms in \eqref{eq:Rhat_expand_2} may have
significant values, which causes the estimates of the signal and
noise subspaces to deviate from the true signal and noise subspaces.

The main idea of our two-step root-MUSIC algorithm is to modify the
sample data covariance matrix at the second step based on the DOA
estimates obtained at the first step. The modified covariance matrix
is obtained by deducting a scaled version of the estimated
undesirable terms from the sample data covariance matrix.

We derive the steps of the proposed method for a general source
covariance matrix $\bS$, so that correlated sources can also be
handled by the algorithm. Furthermore, the proposed method can
also be beneficial in the case that the assumption on no correlation
between the source and noise vectors is not fully met. This is
achieved by estimating and removing the correlation terms between
the source and noise vectors from the sample data covariance matrix.

The steps of the proposed method are listed in Table \ref{alg:1}.
The algorithm starts by computing the sample data covariance matrix
\eqref{eq:Rhat_conventional}. Then, DOAs are estimated using the
root-MUSIC algorithm. The superscript $(\cdot)^{(1)}$ refers to the
estimation made at the first step. At the second step, the
Vandermonde matrix is formed using the available estimates of the
DOAs. Then, the amplitudes of the sources are estimated such that
the squared norm of the differences between the observations and the
estimates are minimized. The corresponding problem is formulated
as
\begin{equation}
\hat{\bs}(t) = \text{arg}~\min_{\bs} \| \bx(t) - \widehat{\bA} \bs
\|_2^2. \label{eq:obj_est_err}
\end{equation}

The minimization of \eqref{eq:obj_est_err} is performed using the
least squares (LS) technique and the corresponding solution is given
as
\begin{equation}
\hat{\bs}(t) = \left(\widehat{\bA}^H \widehat{\bA}\right)^{-1}
\widehat{\bA}^H \bx(t).
\end{equation}
The noise component is then estimated as the difference between the
estimated signal and the observation made by the array, i.e.,
\begin{equation}
\hat{\bn}(t) = \bx(t) - \widehat{\bA} \hat{\bs}(t).
\end{equation}
After estimating the signal and noise vectors, the third term in
\eqref{eq:Rhat_expand_2} can be found as
\begin{eqnarray}
\bT \hspace{-2mm} & \triangleq & \hspace{-2mm} \widehat{\bA} \left\{
\frac{1}{N}\sum_{t=1}^N
\hat{\bs}(t) \hat{\bn}^H(t) \right\} \nonumber \\
& = & \hspace{-2mm} \widehat{\bA} \Bigg\{ \frac{1}{N}\sum_{t=1}^N
\left(\widehat{\bA}^H \widehat{\bA}\right)^{-1} \widehat{\bA}^H
\bx(t) \left( \bx^H(t) - \bx^H(t) \widehat{\bA}
\left(\widehat{\bA}^H \widehat{\bA}\right)^{-1}
\widehat{\bA}^H \right) \Bigg\} \nonumber \\
& = & \hspace{-2mm} \widehat{\bP}_A \left\{ \frac{1}{N}\sum_{t=1}^N
\bx(t) \bx^H(t) \left( \bI_M - \widehat{\bP}_A \right) \right\}
\nonumber \\
& = & \hspace{-2mm} \widehat{\bP}_A \widehat{\bR}
\widehat{\bP}_A^{\bot} \label{eq:T_PARPAbot}
\end{eqnarray}
where
\begin{equation}
\widehat{\bP}_A \triangleq \widehat{\bA} \left(\widehat{\bA}^H
\widehat{\bA}\right)^{-1} \widehat{\bA}^H
\end{equation}
is an estimation for the projection matrix of the signal subspace,
and
\begin{equation}
\widehat{\bP}_A^{\bot} \triangleq \bI_M - \widehat{\bP}_A
\label{eq:PAhatbot_I_PAhat}
\end{equation}
is an estimation for the projection matrix of the noise subspace.
The forth term in \eqref{eq:Rhat_expand_2} is equal to the Hermitian
of the third term, i.e., $\bT^H$. Finally, the modified data
covariance matrix is obtained by deducting a scaled version of the
estimated terms from the initial sample data covariance matrix as
follows
\begin{equation}
\widehat{\bR}^{(2)} = \widehat{\bR} - \gamma \left( \bT+ \bT^H
\right).\label{eq:step2_Rhat2}
\end{equation}

The scaling factor $\gamma$ in \eqref{eq:step2_Rhat2} is a real
number between zero and one. Ideally, the value of $\gamma$ would be
equal to $1$ if the estimates of the undesirable terms were
perfect. However, estimation errors are inevitable, and therefore,
we have introduced $\gamma$ to deal with the imperfections.  The
scaling factor $\gamma$ can be considered as a reliability factor
which takes a value close to $1$ for an estimate of $\bT$ with small
error and a small value if an estimate of $\bT$ is erroneous. Given
a value for $\gamma$, the modified data covariance matrix
$\widehat{\bR}^{(2)}$ is computed and the DOAs are estimated again
using the root-MUSIC algorithm.

The value of $\gamma$ can be fixed to a predetermined value before
running the algorithm, or it can be obtained based on the
observations. Since $\gamma$ is a real number between zero and one,
we can consider different values for $\gamma$ taken on a grid (e.g.
$\gamma = 0,~0.1,~0.2,\cdots,~1$). For each value of $\gamma$, a set
of DOA estimates is obtained based on the modified data covariance
matrix. Next, we determine which value of $\gamma$ results in a
better estimation. This can be done by choosing a set of DOA
estimates that has a higher likelihood of being the set of true
DOAs. In other words, we use the maximum likelihood (ML) criterion
to evaluate the quality of the estimated DOAs. Since the system model
given in \eqref{eq:R_model} is stochastic, we use the stochastic ML
(SML) objective function given by \cite{Stoica90:Stoch_Max_Like}
\begin{equation}
F_{SML}(\gamma) = \ln \det \left( \widehat{\bP}_A^{(2)}
\widehat{\bR} \widehat{\bP}_A^{(2)} + \frac{\text{Tr}
\left\{\widehat{\bP}_A^{\bot(2)} \widehat{\bR} \right\} }{M - K}
\widehat{\bP}_A^{\bot(2)} \right) \label{eq:SML_obj_func}
\end{equation}
where $\text{Tr} \left\{ \cdot \right\}$ stands for the trace
operator, $\widehat{\bP}_A^{(2)}$ is an estimation of the projection
matrix of the signal subspace obtained from the estimated DOAs based
on the modified data covariance matrix and
$\widehat{\bP}_A^{\bot(2)} = \bI_M - \widehat{\bP}_A^{(2)}$. The
objective function in \eqref{eq:SML_obj_func} is evaluated for each
value of $\gamma$. Then, the set of DOA estimates corresponding to
the value of $\gamma$ that minimizes \eqref{eq:SML_obj_func} is
chosen as the output of the algorithm.

\begin{table}[t]%\vspace{2mm}
\caption{Two-step root-MUSIC algorithm}\label{alg:1}\vspace{-5mm}
\begin{center} \normalsize{
\begin{tabular}{l}
\hline
%Iterative root-MUSIC Algorithm\\
%\hline
\textbf{Inputs}:\\
\hspace{1mm} $M,~d,~\lambda,~N,~K,$ and\\
\hspace{1mm} received vectors $\bx(1),~\bx(2),\cdots,~\bx(N)$\\
\textbf{Outputs}:\\
\hspace{1mm} Estimates $\hat{\theta}_1^{(2)},~\hat{\theta}_2^{(2)},
\cdots,~\hat{\theta}_K^{(2)}$ \vspace{0.5mm} \\
\hline
\textbf{Step 1}:\\
\hspace{1mm} $\widehat{\bR} =
\frac{1}{N}\sum_{t=1}^N \bx(t) \bx^H(t)$ \vspace{0.5mm} \\
\hspace{1mm}
$\left\{\hat{\theta}_1^{(1)},~\hat{\theta}_2^{(1)},\cdots,~\hat{\theta}_K^{(1)}
\right\}
\leftarrow \text{root-MUSIC}\left(\widehat{\bR},K,d,\lambda \right)$\\
\textbf{Step 2:}\\
\hspace{1mm} $\widehat{\bA} =
\left[\ba\left(\hat{\theta}_1^{(1)}\right),~
\ba\left(\hat{\theta}_2^{(1)}\right),\cdots,~
\ba\left(\hat{\theta}_K^{(1)}\right)\right]$ \vspace{1.0mm} \\
\hspace{1mm} $\widehat{\bP}_A = \widehat{\bA} \left(\widehat{\bA}^H
\widehat{\bA}\right)^{-1} \widehat{\bA}^H$ \\
\hspace{1mm} $\widehat{\bP}_A^{\bot} = \bI_M -
\widehat{\bP}_A$ \\
\hspace{1mm} $\bT = \widehat{\bP}_A
\widehat{\bR} \widehat{\bP}_A^{\bot} $\\
\hspace{1mm} \textbf{Determine $\gamma$ as the minimizer of}
\eqref{eq:SML_obj_func} \\
\hspace{2mm} $\widehat{\bR}^{(2)} = \widehat{\bR} - \gamma \left(
\bT+
\bT^H \right)$\\
\hspace{2mm}
$\left\{\hat{\theta}_1^{(2)},~\hat{\theta}_2^{(2)},
\cdots,~\hat{\theta}_K^{(2)}\right\}
\leftarrow \text{root-MUSIC}\left(\widehat{\bR}^{(2)},
K, d, \lambda \right)$ \vspace{1mm}\\
\hline
\end{tabular}}
\end{center}
\end{table}

\subsection{Root-swap root-MUSIC algorithm}
Consider the root-MUSIC polynomial $\ba^T(z^{-1}) \bG \bG^H \ba(z)$
which is formed by the noise eigenvectors obtained from the
eigendecomposition of the data covariance matrix $\bR$. This
polynomial has $K$ number of roots on the unit circle which
correspond to the signal sources. Let these $K$ roots be denoted by
$z_1,~z_2,\cdots,~z_K$ and be referred to as the true signal roots.
The polynomial also has additional $M - K - 1$ number of roots
inside the unit circle. Let these roots be referred to as the true
noise roots and be denoted by $z_{K+1},~z_{K+2},\cdots,~z_{M-1}$.

An estimation for the root-MUSIC polynomial can be formed using the
noise eigenvectors obtained from the sample data covariance matrix
$\widehat{\bR}$. Let us assume that in the estimation of the noise
and signal subspaces, no subspace swap has occurred
\cite{Scharf95:prob_subspace_swap}. The estimated polynomial is
given by $\ba^T(z^{-1}) \widehat{\bG} \widehat{\bG}^H \ba(z)$. This
polynomial has $M-1$ number of roots inside the unit circle. Let
$\hat{z}_1,~\hat{z}_2,\cdots,~\hat{z}_{K}$ be the roots of the
estimated root-MUSIC polynomial which correspond to
$z_1,~z_2,\cdots,~z_K$. We refer to these roots as the estimated
signal roots. Furthermore, let
$\hat{z}_{K+1},~\hat{z}_{K+2},\cdots,~\hat{z}_{M-1}$ be the roots
corresponding to $z_{K+1},~z_{K+2},\cdots,~z_{M-1}$. These roots are
referred to as the estimated noise roots.

In the root-MUSIC method, we do not have the knowledge about which
of the roots of the estimated root-MUSIC polynomial correspond to
the true signal roots. The conventional rule is to select $K$ number
of the estimated roots which are closer to the unit circle as the
estimates for the true signal roots. Then, the DOAs are estimated
based on the angles of these roots.

Due to the finiteness of the available samples, the estimated roots
obtained from the sample data covariance matrix $\widehat{\bR}$
deviate from their corresponding true roots obtained from the true
data covariance matrix $\bR$. Let $r_i$ and $\hat{r}_i$ represent
the magnitudes of $z_i$ and $\hat{z}_i$ for $1 \leq i \leq M-1$,
respectively. Furthermore, let $\Delta r_i \triangleq \hat{r}_i -
r_i$ be the difference between the magnitude of the $i$-th estimated
root and the magnitude of the corresponding true root. It is shown
in \cite{Proakis92:root_Music_performance} that $\Delta r_i$ (for
the signal roots) has a variance which is proportional to
$\sigma_{\text{n}}^2 / N$. Therefore, $\Delta r_i$ can have a
significant value for a small number of samples and a large value of
$\sigma_{\text{n}}^2$ (low SNR region). Consequently, there can be a
considerable probability that an estimated signal root takes a
smaller magnitude than an estimated noise root. We refer to this
phenomenon as a root-swap. The root-swap probability is
approximately found in Appendix~\ref{sec:appndx_root_swap_prob} as
\begin{eqnarray}
&& \hspace{-15mm} P(\text{root-swap}) \approx 1 - \prod_{k=1}^K
\prod_{m=K+1}^{M-1} Q\left( \frac{-1 + r_m + \sigma_k
\sqrt{M-K-(3/4)} }{\sqrt{\sigma_k^2 / 4}} \right)
\label{eq:root_swap_prob_approx_final}
\end{eqnarray}
where $Q\left( \cdot \right)$ is the tail probability of the
standard normal distribution and $\sigma_k^2/4$ is the variance of
$\Delta r_k$, and it is proportional to $\sigma_{\text{n}}^2 / N$.

In the case that the root-swap happens, selecting the first closest
$K$ roots to the unit circle results in picking a noise root instead
of a signal root. To deal with this problem, we propose an algorithm
that considers different combinations of the roots as candidates for
signal roots. The method is dubbed the root-swap root-MUSIC
algorithm.

The root-MUSIC polynomial has $M-1$ number of roots inside the unit
circle. Our goal is to find the roots which have a higher likelihood
of being associated with the $K$ sources. Consider choosing $K$
number of roots out of the $M-1$ roots inside the unit circle. There
are $N_c \triangleq (M - 1)! / \left( K! (M - K - 1)! \right)$
different possible combinations. Let $\Gamma \triangleq \left\{
\Theta_1,~\Theta_2,~\cdots,~\Theta_{N_c} \right\}$ where $\Theta_i$
($1 \leq i \leq N_c$) is a set containing the DOA estimates obtained
from the $i$-th combination of the roots. Then, the root-swap
root-MUSIC method estimates the DOAs as
\begin{equation}
\left\{ \hat{\theta}_1,~\hat{\theta}_2,\cdots,~\hat{\theta}_K
\right\} = \text{arg}~\min_{\Theta \in \Gamma} F_{SML}\left( \Theta
\right)
\end{equation}
where $F_{SML}\left( \Theta \right)$ is the SML function given by
\begin{equation}
F_{SML}(\Theta) = \ln \det \left( \bP_{\Theta} \widehat{\bR}
\bP_{\Theta} + \frac{\text{Tr} \left\{\bP_{\Theta}^{\bot}
\widehat{\bR} \right\} }{M - K} \bP_{\Theta}^{\bot} \right)
\end{equation}
and $\bP_{\Theta}$ is the signal projection matrix obtained from
$\Theta$ as
\begin{equation}
\bP_{\Theta} \triangleq \bA(\Theta) \left(\bA^H(\Theta)
\bA(\Theta)\right)^{-1} \bA^H(\Theta).
\end{equation}

The complexity of the introduced root-swap root-MUSIC method can be
reduced by pre-eliminating some of the roots. Specifically, let $p
\leq K$ roots closest to the unit circle be picked, and let $q$
number of roots closest to the origin (furthest from the unit
circle) be ignored. Our task is to choose $K - p$ number of roots
out of $M - p - q -1$ roots. Then, there are $N_{\rm r} \triangleq
(M - p - q - 1)! / \left( (K - p)! (M - K - q - 1)! \right)$
different possible combinations which is significantly smaller
than $N_{\rm c}$. The rest of the algorithm is the same as above
except for that here each combination contains $K - p$ number
of roots. Therefore, in order to evaluate the SML function, the
fixed $p$ pre-selected roots are added to each combination.

\section{Subspace Leakage}
\label{sec:Sub_leak}
The performance breakdown of the subspace based DOA estimation
methods in the threshold region has been associated with the
subspace leakage. In this section, we study the subspace leakage
for both steps of the proposed two-step root-MUSIC algorithm.

\subsection{Definition}
Recall the matrices $\bG$ and $\bE$ which are composed of the true
noise and signal eigenvectors obtained from the eigendecomposition
of the data covariance matrix $\bR$. Note that the matrix of the
eigenvectors $\bQ_R = \left[ \bG~\bE \right] \in \mathbb{C}^{M
\times M}$ is a unitary matrix $\left( \bQ_R \bQ_R^H = \bI_M
\right)$, therefore
\begin{equation}
\bG \bG^H + \bE \bE^H = \bI_M
\end{equation}
or
\begin{equation}
\bP^{\bot} + \bP = \bI_M \label{eq:P_orth+P=I}
\end{equation}
where, $\bP^{\bot} \triangleq \bG \bG^H$ and $\bP \triangleq \bE
\bE^H$ are the true projection matrices of the noise and signal
subspaces.

Ideally, the estimation of each signal eigenvector $\hat{\be}_k$ $(1
\leq k \leq K)$ would perfectly fall in the true signal subspace. In
practice, however, the energy of the projection of $\hat{\be}_k$
into the noise subspace $\| \bP^{\bot} \hat{\be}_k \|_2^2$ is almost
surely nonzero, which can be viewed as the leakage of $\hat{\be}_k$
into the true noise subspace.

We define the subspace leakage as the average value of the energy of
the estimated signal eigenvectors leaked into the true noise
subspace, i.e.,
\begin{equation}
\rho \triangleq \frac{1}{K} \sum_{k = 1}^K \| \bP^{\bot} \hat{\be}_k
\|_2^2.
\end{equation}
Note that $\bP^{\bot}$ is the orthogonal projection matrix.
Therefore, $\rho$ can be written as
\begin{equation}
\rho = \frac{1}{K} \sum_{k = 1}^K \hat{\be}_k^H \bP^{\bot}
\hat{\be}_k. \label{eq:sim1_rho_def}
\end{equation}
Using \eqref{eq:P_orth+P=I} and some algebra, the expression
\eqref{eq:sim1_rho_def} can be simplified to
\begin{eqnarray}
\rho \hspace{-2mm} & = & \hspace{-2mm} \frac{1}{K} \sum_{k = 1}^K
\hat{\be}_k^H \left( \bI_M - \bP \right)
\hat{\be}_k \nonumber \\
& = & \hspace{-2mm} 1 - \frac{1}{K} \sum_{k = 1}^K \text{Tr} \left\{
\hat{\be}_k \hat{\be}_k^H
\bP \right\} \nonumber \\
& = & \hspace{-2mm} 1 - \frac{1}{K} \text{Tr} \left\{ \left( \sum_{k
= 1}^K \hat{\be}_k \hat{\be}_k^H \right)
\bP \right\} \nonumber \\
& = & \hspace{-2mm} 1 - \frac{1}{K} \text{Tr} \left\{ \widehat{\bP}
\bP \right\} \label{eq:sub_leak_def}
\end{eqnarray}
where $\widehat{\bP} \triangleq \widehat{\bE} \widehat{\bE}^H$ is
the estimated signal projection matrix.

\subsection{Analysis of two-step root-MUSIC algorithm}
The estimated signal and noise projection matrices obtained from the
eigendecomposition of the sample data covariance matrix
$\widehat{\bR}$ are deviated from the true signal and noise
projection matrices. Let $\Delta \bR \triangleq \widehat{\bR} - \bR$
be the estimation error of the data covariance matrix, and let
\begin{eqnarray}
\bV \hspace{-2mm} & \triangleq & \hspace{-2mm}
\bR - \sigma_{\text{n}}^2\bI_M = \bA\bS\bA^H \nonumber \\
& = & \hspace{-2mm} \sum_{k = 1}^K \left( \lambda_{M - K + k} -
\sigma_{\text{n}}^2 \right) \be_k \be_k^H. \label{eq:V_R_sigma2I}
\end{eqnarray}
Denote the pseudo-inverse of $\bV$ as $\bV^{\dag} \in \mathbb{C}^{M
\times M}$. It is given by
\begin{equation}
\bV^{\dag} = \sum_{k = 1}^K \frac{1}{\lambda_{M - K + k} -
\sigma_{\text{n}}^2} \be_k \be_k^H. \label{eq:pseudoinverse_V}
\end{equation}

Let $\rho_1$ and $\rho_2$ be the subspace leakage due to the error
in the estimation of the signal and noise subspaces obtained from
$\widehat{\bR}$ and $\widehat{\bR}^{(2)}$, respectively. Note that
$\rho_1$ only depends on $\bR$ and $\Delta \bR$, and it is not
specific to the proposed two-step root-MUSIC algorithm.

It is shown in Appendix~\ref{sec:appndx_SL1} that $\rho_1$ and its
expected value are given by
\begin{equation}
\rho_1 = \frac{1}{K} \text{Tr} \left\{ \bV^{\dag} \Delta \bR
\bP^{\bot} \Delta \bR \bV^{\dag} \right\}
\end{equation}
and
\begin{equation}
E\left\{ \rho_1 \right\} = \frac{\sigma_{\text{n}}^2 \left( M - K
\right)}{NK} \sum_{k = 1}^K \frac{\lambda_{M - K + k}}{\left(
\lambda_{M - K + k} - \sigma_{\text{n}}^2 \right)^2 }
\label{eq:E_rho1_final_simplified}
\end{equation}
respectively.

It can be seen from \eqref{eq:E_rho1_final_simplified} that the
expected value of the subspace leakage is proportional to
$\sigma_{\text{n}}^2 /N$. Therefore, the amount of the subspace
leakage can be significant for a small number of samples or low SNR
values. The variance of $\rho_1$ has also been studied in
\cite{MSH15:Sub_Leak}, and it has been shown that $\text{Var}
\left(\rho_1 \right)$ is in the order of $1/N^2$.

The subspace leakage at the second step of the two-step root-MUSIC
algorithm is computed in Appendix~\ref{sec:appndx_SL2} and is given
by
\begin{eqnarray}
\rho_2 \hspace{-2mm} & = & \hspace{-2mm} \left( 1 - 2\gamma +
\gamma^2 \right) \rho_1 + \frac{ 2 \left( \gamma - \gamma^2 \right)
}{K} Re \left\{ \text{Tr} \left\{ \bV^{\dag} \Delta \bR \bP^{\bot}
d\bP \right\} \right\} + \frac{\gamma^2}{K} \text{Tr} \left\{ d\bP
\bP^{\bot} d\bP \right\}
\end{eqnarray}
where $Re \left\{ \cdot \right\}$ stands for the real part operator,
and $d\bP$ is the first order term in the Taylor series expansion of
$\widehat{\bP}_A$ around the true DOAs. It is also shown in
Appendix~\ref{sec:appndx_SL2} that the expected value of $\rho_2$
for a fixed value of $\gamma$ is given by
\begin{eqnarray}
E \left\{ \rho_2 \right\} \hspace{-2mm} & = & \hspace{-2mm} \left( 1
- 2\gamma + \gamma^2 \right) E \left\{ \rho_1 \right\} \nonumber \\
&& \hspace{-2mm} + \frac{ 2 \left( \gamma - \gamma^2 \right)
\sigma_{\text{n}}^2 }{NK} Re \left\{ \sum_{k = 1}^K
\frac{\ba_k^{(1)H} \bP^{\bot} \frac{\partial \bA}{\partial \omega_k}
\left( \bA^H \bA \right)^{-1} \bA^H \bV^{\dag} \bR \bV^{\dag}
\ba_k}{2 j \left( \ba_k^{(1)H} \bP^{\bot} \ba_k^{(1)} \right)}
\right\} \nonumber \\
&& \hspace{-2mm} + \frac{\gamma^2 \sigma_{\text{n}}^2}{2NK} \sum_{k
= 1}^K \sum_{i = 1}^K \frac{\text{Tr} \left\{ \left( \frac{\partial
\bA}{\partial \omega_k} \right)^H \bP^{\bot} \frac{\partial
\bA}{\partial \omega_i} \left( \bA^H \bA \right)^{-1} \right\}}{
\left( \ba_k^{(1)H} \bP^{\bot} \ba_k^{(1)} \right) \left(
\ba_i^{(1)H} \bP^{\bot} \ba_i^{(1)} \right)} Re \left\{ \ba_i^H
\bV^{\dag} \bR \bV^{\dag} \ba_k \ba_k^{(1)H} \bP^{\bot} \ba_i^{(1)}
\right\} \nonumber \\
\label{eq:E_rho2_final_simplified}
\end{eqnarray}
where $\omega_k \triangleq 2 \pi (d/\lambda) \sin(\theta_k)$,
$\ba_k$ is a shorthand notation for $\ba(\theta_k)$, and
$\ba_k^{(1)} \in \mathbb{C}^{M \times 1}$ is defined as
\begin{equation}
\ba_k^{(1)} \triangleq \hspace{-1mm} - \hspace{-1mm} \left[0,~e^{-j
\omega_k},~2 e^{-j 2 \omega_k},~\cdots,~(M-1) e^{-j (M-1) \omega_k}
\right]^T. \label{eq:a1_k}
\end{equation}

It can be seen in \eqref{eq:E_rho2_final_simplified} that for
$\gamma = 0$, $E \left\{ \rho_2 \right\}$ reduces to $E \left\{
\rho_1 \right\}$ as expected, and for $\gamma = 1$, the first two
terms in \eqref{eq:E_rho2_final_simplified} are equal to zero.

\section{Numerical Examples and Simulation Results}
\label{sec:Sim_Sub_Leak}
In this section, the performance of the proposed two-step
root-MUSIC and the root-swap root-MUSIC algorithms is investigated
and compared with the performance of the unitary root-MUSIC method
\cite{Gershman00:unitary_root_Music} and the improved unitary
root-MUSIC algorithm based on pseudo-noise resampling
\cite{Qian14:impro_uni_r_Music_Pseudo_noise}. We also consider the
combination of the proposed methods with the other methods in order
to achieve further performance improvement. Compared to the
root-MUSIC method, the unitary root-MUSIC algorithm has a lower
computational complexity as it uses the eigendecomposition of a
real-valued covariance matrix. Furthermore, the unitary root-MUSIC
algorithm has better performance for the case that the sources are
correlated. The improved unitary root-MUSIC algorithm based on
pseudo-noise resampling increases the estimator complexity, but it
is advantageous in removing the outliers, which results in better
performance.

We consider $K = 2$ sources impinging on an array of $M = 10$
antenna elements from directions $\theta_1 = 35\,^{\circ} \times
(\pi/180)$ and $\theta_2 = 37\,^{\circ} \times (\pi/180)$. The
interelement spacing is set to $d = \lambda/2$ and the number of
snapshots is $N = 10$. Each source vector $\bs(t)$ is considered to
be independent from the source vectors at other time instances and
to have the circularly-symmetric complex jointly-Gaussian
distribution $\mathcal{N}_C(0,\bS)$. The source covariance matrix
$\bS$ is given by
\begin{equation}
\bS = \sigma_{\text{s}}^2 \left[ \begin{array}{cc}
1 & r \\
r & 1
\end{array} \right] \nonumber
\end{equation}
where $0 \leq r \leq 1$ is the correlation coefficient. The SNR is
defined as $\text{SNR} \triangleq
10\log_{10}\left(\sigma_{\text{s}}^2/\sigma_{\text{n}}^2\right)$.

The performance of the proposed algorithms is investigated by
considering the subspace leakage, mean squared error (MSE),
probability of source resolution, and conditional mean squared error
(CMSE). Source resolution is defined as the event when both DOAs are
estimated within one degree of their corresponding true values,
i.e., the difference between the true value of each DOA and its
estimated value is less than $1\,^{\circ} \times (\pi/180)$. The
CMSE is defined as the expected value of the estimation error
conditioned on successful source resolution, i.e.,
\raisebox{0pt}[5mm][3mm]{$E\left\{\sum_{k=1}^K\| \hat{\theta}_k -
\theta_k \|_2^2 ~\Big|~ \text{successful source resolution}
\right\}$}. The reason for using the CMSE is to further investigate
the accuracy of the algorithms after making successful detection. We
estimate the probability of root-swap, subspace leakage, MSE,
probability of source resolution, and CMSE using the Monte Carlo
method with $10^5$ number of trials. Two cases are considered in the
simulations: 1) the two sources are uncorrelated, i.e., $r = 0$, and
2) the two sources are correlated with a correlation coefficient of
$r = 0.9$.

\begin{figure}[t]
\psfrag{SNR}{SNR (dB)} \psfrag{RS}{Probability}
    \includegraphics[width=25cm]{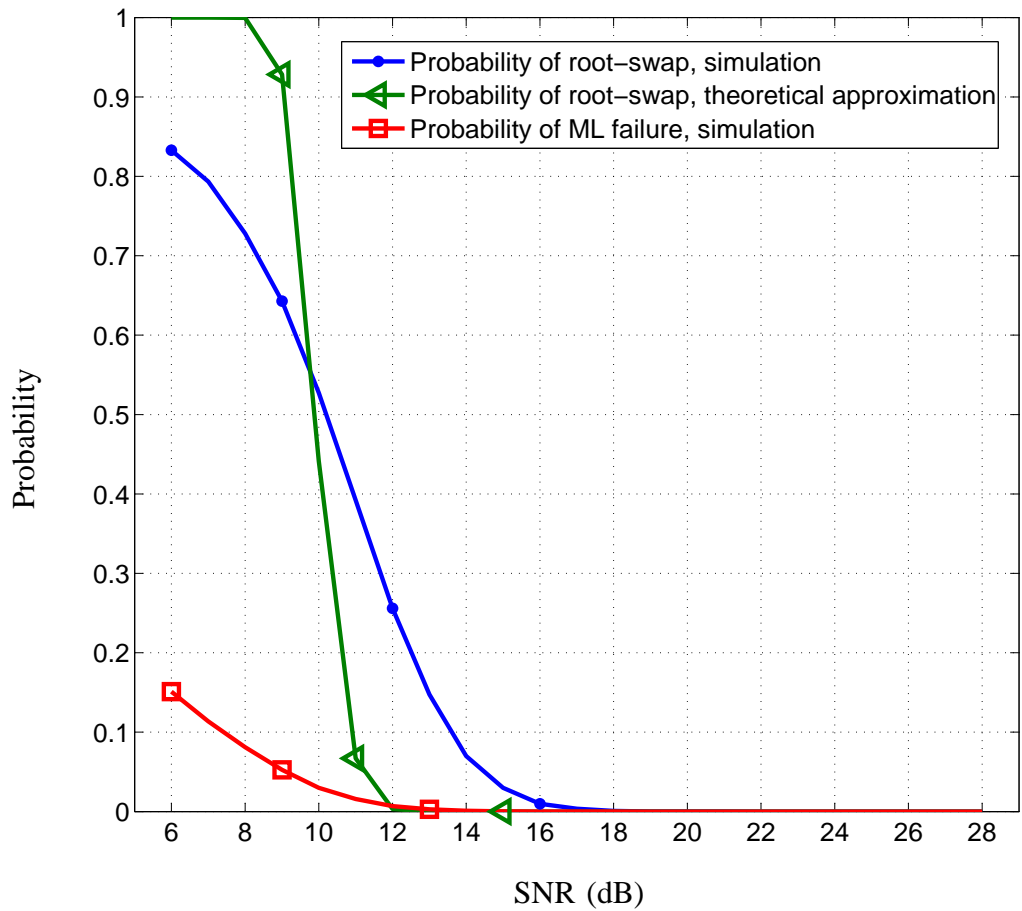}
\vspace{-5mm} \caption{Probability of root-swap and probability of
ML failure versus SNR for uncorrelated sources.
\label{fig:RS_M10_N10_c0}}\vspace{-3mm}
\end{figure}

Let us start by investigating the probability of root-swap in the
root-MUSIC algorithm for the case of the uncorrelated sources. The
probability of root-swap is estimated using the Monte Carlo
simulations. Its approximate value is also obtained using
\eqref{eq:root_swap_prob_approx_final}. The corresponding curves are
shown in Fig.~\ref{fig:RS_M10_N10_c0}. It can be seen that at the
low SNR region, the chance that a root-swap occurs is quite
significant, which results in the performance breakdown of the
root-MUSIC algorithm. This problem justifies the need for a method
to deal with the root-swap phenomenon. In this paper, we proposed
the root-swap root-MUSIC algorithm which instead of picking the
roots closer to the unit circle, selects the roots based on the SML
criterion. In Fig.~\ref{fig:RS_M10_N10_c0}, we thus also draw a
curve which shows the probability that the selected roots by the ML
criterion include a noise root. This situation is considered as a
breakdown, and therefore, the corresponding probability is called
the probability of ML failure. As can be seen, this probability is
significantly smaller than the probability of root-swap. As a
result, it is expected that the root-swap root-MUSIC algorithm
outperforms the conventional root-MUSIC method. This will be shown
in the rest of this section.

\begin{figure}[t]
\psfrag{SNR}{SNR (dB)} \psfrag{SL}{Subspace leakage (dB)}
    \includegraphics[width=25cm]{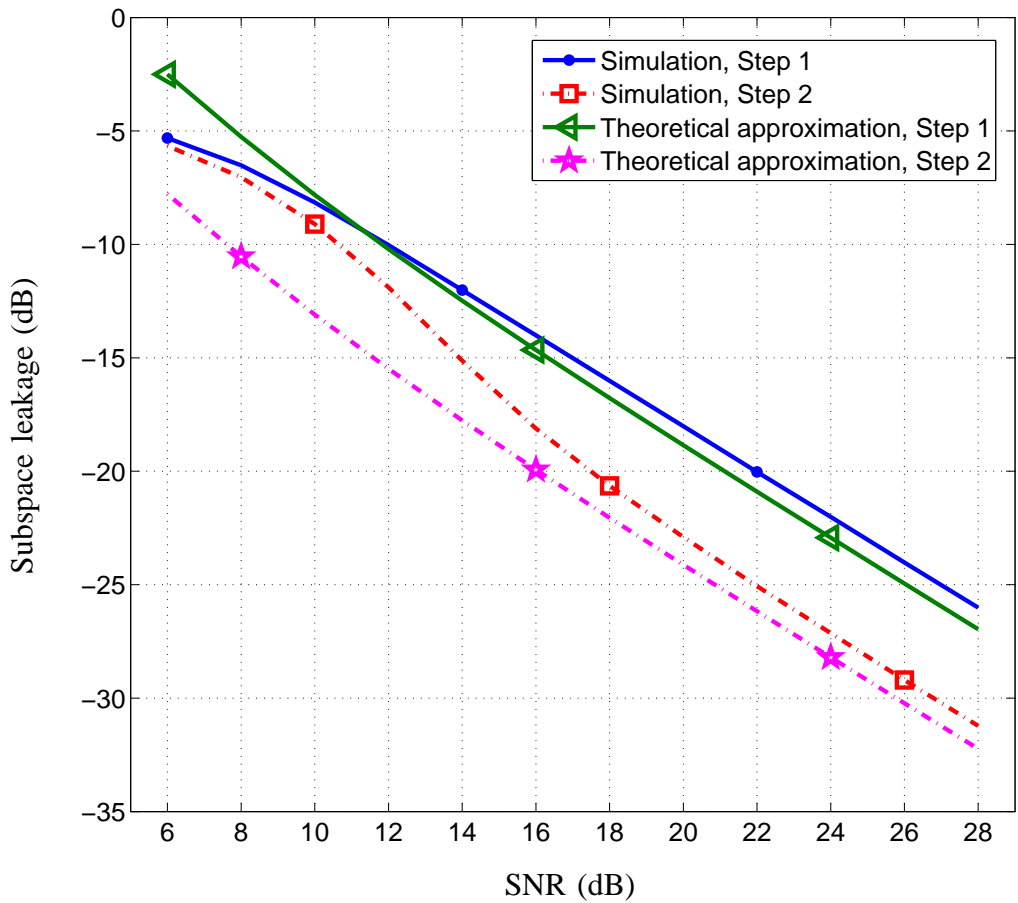}
\vspace{-5mm} \caption{Subspace leakage versus SNR for uncorrelated
sources. The solid and dashed lines represent the subspace leakage
at the first and second steps of the proposed two-step root-MUSIC
algorithm, respectively. \label{fig:SL_M10_N10_c0}}\vspace{-3mm}
\end{figure}

The subspace leakage in the two-step root-MUSIC algorithm for the
case of the uncorrelated sources is investigated next. The expected
value of the subspace leakage is estimated using
\eqref{eq:sub_leak_def} and the Monte Carlo simulations. The
approximate value for the subspace leakage is also obtained from the
theoretical derivations in \eqref{eq:E_rho1_final_simplified} and
\eqref{eq:E_rho2_final_simplified}. The value of $\gamma$ is fixed
at $0.5$. The results are shown in Fig.~\ref{fig:SL_M10_N10_c0}. The
solid lines represent the subspace leakage at the first step, and
the dashed lines depict the subspace leakage at the second step of
the proposed two-step root-MUSIC algorithm. It can be seen that the
curves obtained from the simulations are very close to those
obtained from our theoretical derivations at high SNR values. At the
low SNR region, the curve associated with the theoretical
approximation at the second step deviates from the curve obtained by
simulations. The reason is that in the derivations, the first order
Taylor series expansion is used. More accurate results can be
obtained by using higher order Taylor series. However, the
computations can become intractable. In
Fig.~\ref{fig:SL_M10_N10_c0}, it can be observed from both
theoretical and simulation results that the subspace leakage from
the modified covariance matrix at the second step is significantly
smaller than the subspace leakage from the sample data covariance
matrix at the first step. This is achieved by removing the
undesirable terms from the sample data covariance matrix leading to
an estimate of the signal projection matrix that is closer to the
true signal projection matrix, which is equivalent to a lower
subspace leakage at the second step.

We next consider the performance of the proposed two-step algorithm
when applied to the root-MUSIC \cite{Barabell83:root_Music}, unitary
root-MUSIC \cite{Gershman00:unitary_root_Music}, improved unitary
root-MUSIC with pseudo-noise resampling
\cite{Qian14:impro_uni_r_Music_Pseudo_noise}, root-swap unitary
root-MUSIC, and root-swap unitary root-MUSIC with pseudo-noise
resampling methods. The unitary root-MUSIC algorithm takes benefit
from the forward-backward averaging \cite{Pillai89:Forward_Backward}
which is approximately equivalent to doubling the number of samples.
For the cases that the pseudo-noise resampling is used, $P$
represents the number of times that the resampling process has been
performed. In the figures, the root-MUSIC, unitary root-MUSIC, and
root-swap unitary root-MUSIC methods are denoted by R-MUSIC,
UR-MUSIC, and RSUR-MUSIC, respectively. The value of the scaling
factor $\gamma$ is obtained by minimizing the SML function as
described in the two-step root-MUSIC method. In the root-swap
algorithm, the parameters $p$ and $q$ are set to $p = 1$ and $q =
0$, which means the closest root to the unit circle is picked up and
paired with other roots one at a time in order to find the pair of
DOA estimates that minimizes the SML function. In this case, the
number of different possible combinations of the roots is $N_r = 8$.
The number of samples used for the pseudo-noise resampling method is
set to $P = 50$. According to our simulations, using more number of
samples would not yield in any considerable improvement in the
performance.

\begin{figure}[t]
\psfrag{SNR}{SNR (dB)} \psfrag{MSE}{MSE (dB)}
    \includegraphics[width=25cm]{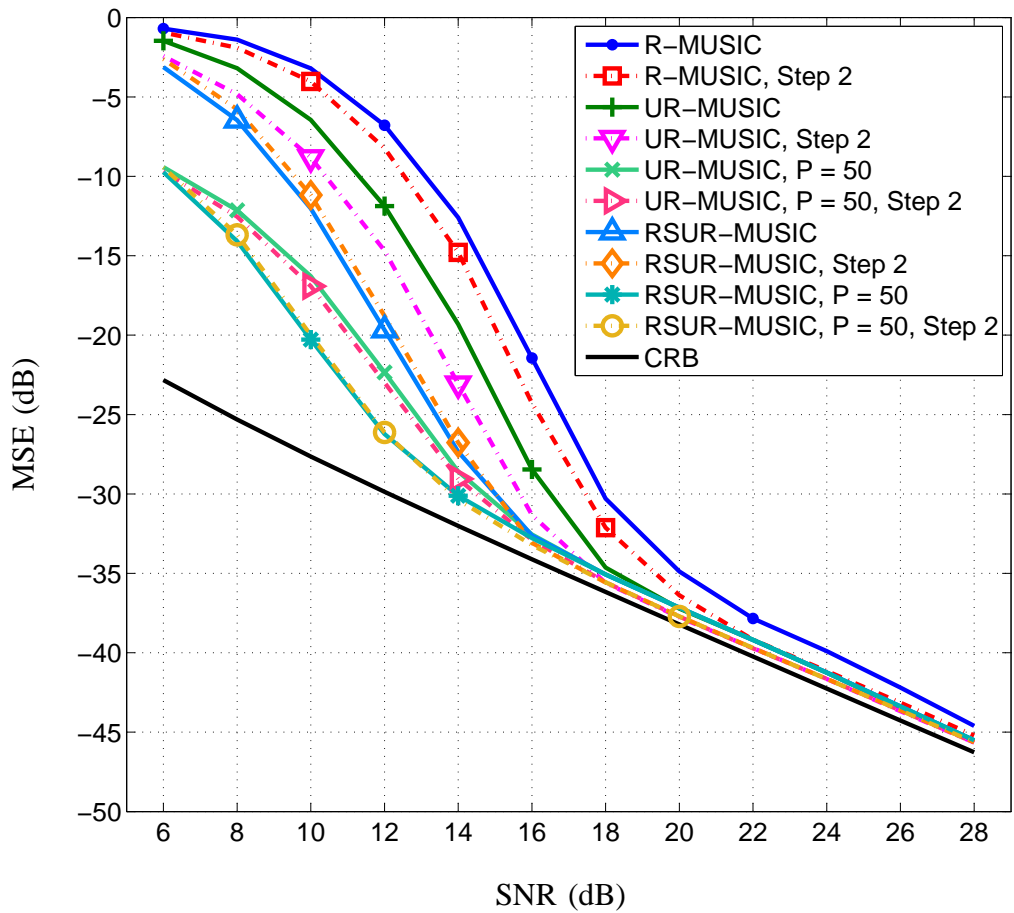}
\vspace{-5mm} \caption{MSE versus SNR for uncorrelated sources. The
solid and dashed lines are based on the first and second steps of
the proposed two-step method, respectively. The methods used in the
two-step algorithm are R-MUSIC, UR-MUSIC, and RSUR-MUSIC methods.
$P$ is the number of samples used for the pseudo-noise resampling
algorithm. \label{fig:MSE_M10_N10_c0}}\vspace{-3mm}
\end{figure}

The MSE versus SNR performance of the methods tested for the case of
the uncorrelated sources is presented in
Fig.~\ref{fig:MSE_M10_N10_c0}. The corresponding CRB
\cite{Stoica01:Stocastic_CRB_array_proc} is also shown in the
figure. For the R-MUSIC method, the modification of the covariance
matrix in the second step of the introduced two-step method shifts
the MSE curve by almost half~a~dB to the left. For the UR-MUSIC
method the improvement is more significant and is about one~dB. For
the rest of the methods, there is no considerable change in the MSE
performance. However, as it will be shown in the next figures, the
modification of the covariance matrix has benefits in terms of the
CMSE performance and probability of source resolution for these
methods. It can also be seen from Fig.~\ref{fig:MSE_M10_N10_c0} that
the proposed RSUR-MUSIC algorithm performs about $2$~dB better than
the UR-MUSIC method, while imposing only a small amount of
computational complexity for evaluating the SML function for $N_r =
8$ different combinations of the roots. The best performance is
achieved by the RSUR-MUSIC algorithm combined with the pseudo-noise
resampling method.

\begin{figure}[t]
\psfrag{SNR}{SNR (dB)} \psfrag{PD}{Probability of source resolution}
    \includegraphics[width=25cm]{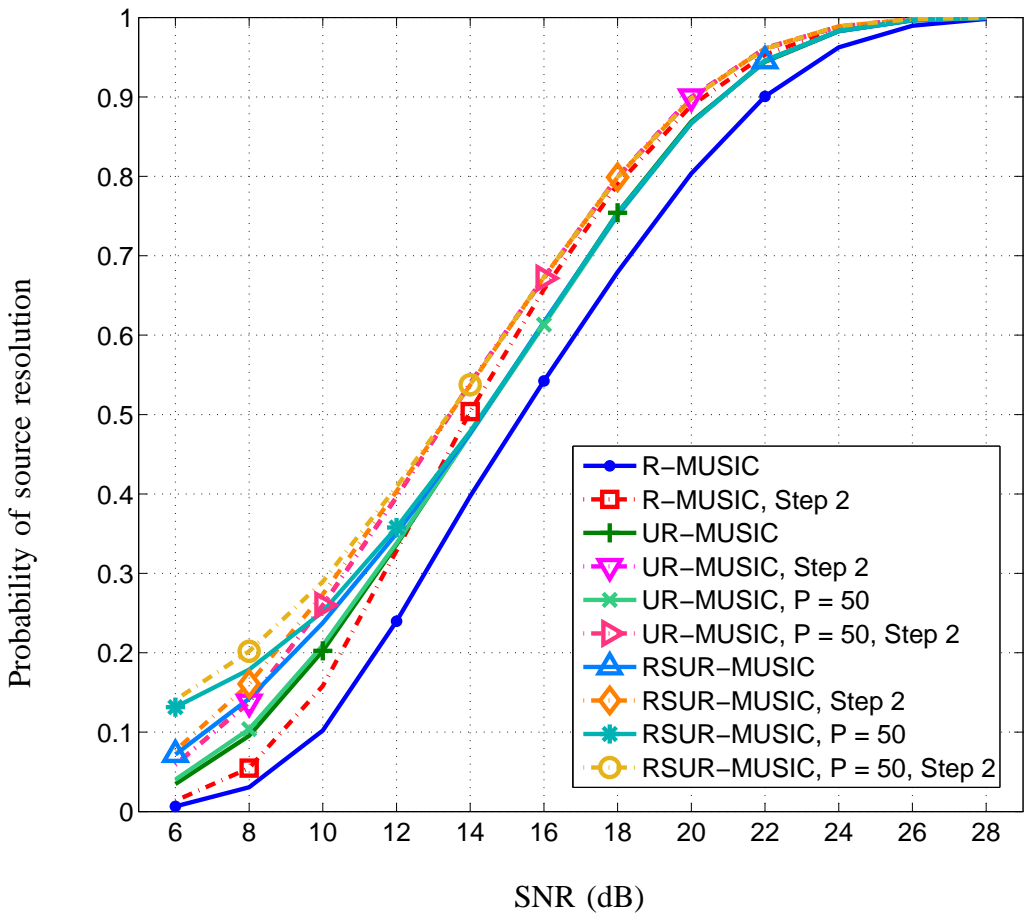}
\vspace{-5mm} \caption{Probability of source resolution versus SNR
for uncorrelated sources. The solid and dashed lines are based on
the first and second steps of the proposed two-step method,
respectively. The methods used in the two-step algorithm are
R-MUSIC, UR-MUSIC, and RSUR-MUSIC methods.
\label{fig:PD_M10_N10_c0}}\vspace{-3mm}
\end{figure}

Fig.~\ref{fig:PD_M10_N10_c0} shows probability of source resolution
versus SNR for the uncorrelated sources. For the R-MUSIC method, the
second step of the two-step algorithm improves the performance by
$1$ to $2$~dB. The rest of the algorithms have almost the same
performance with the root-swap based methods slightly outperforming
the other algorithms at low SNR values. It is observed that the
second step of the two-step algorithm results in about $1$~dB
improvement in the performance.

\begin{figure}[t]
\psfrag{SNR}{SNR (dB)} \psfrag{CMSE}{CMSE (dB)}
    \includegraphics[width=25cm]{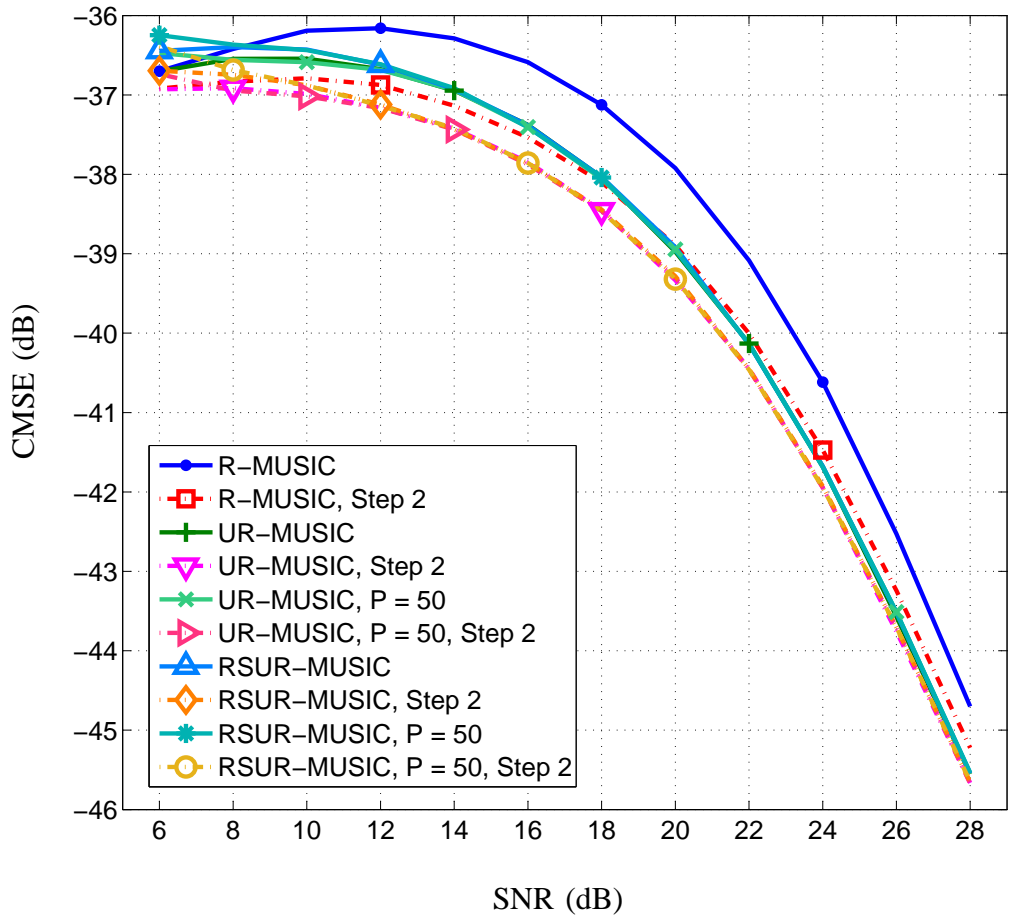}
\vspace{-5mm} \caption{CMSE versus SNR for uncorrelated sources. The
solid and dashed lines are based on the first and second steps of
the proposed two-step method, respectively. The methods used in the
two-step algorithm are R-MUSIC, UR-MUSIC, and RSUR-MUSIC methods.
\label{fig:CMSE_M10_N10_c0}}\vspace{-3mm}
\end{figure}

Finally, Fig.~\ref{fig:CMSE_M10_N10_c0} illustrates the performance
of the algorithms tested for the uncorrelated sources in terms of
the CMSE. The R-MUSIC method is significantly improved by the
two-step method with an improvement ranging from $5$~dB at low SNR
values to $1$~dB at high SNR values. The rest of the algorithms show
similar performance, and the application of the two-step method
leads to up to $2$~dB improvement in the CMSE performance.

\begin{figure}[t]
\psfrag{SNR}{SNR (dB)} \psfrag{RS}{Probability}
    \includegraphics[width=25cm]{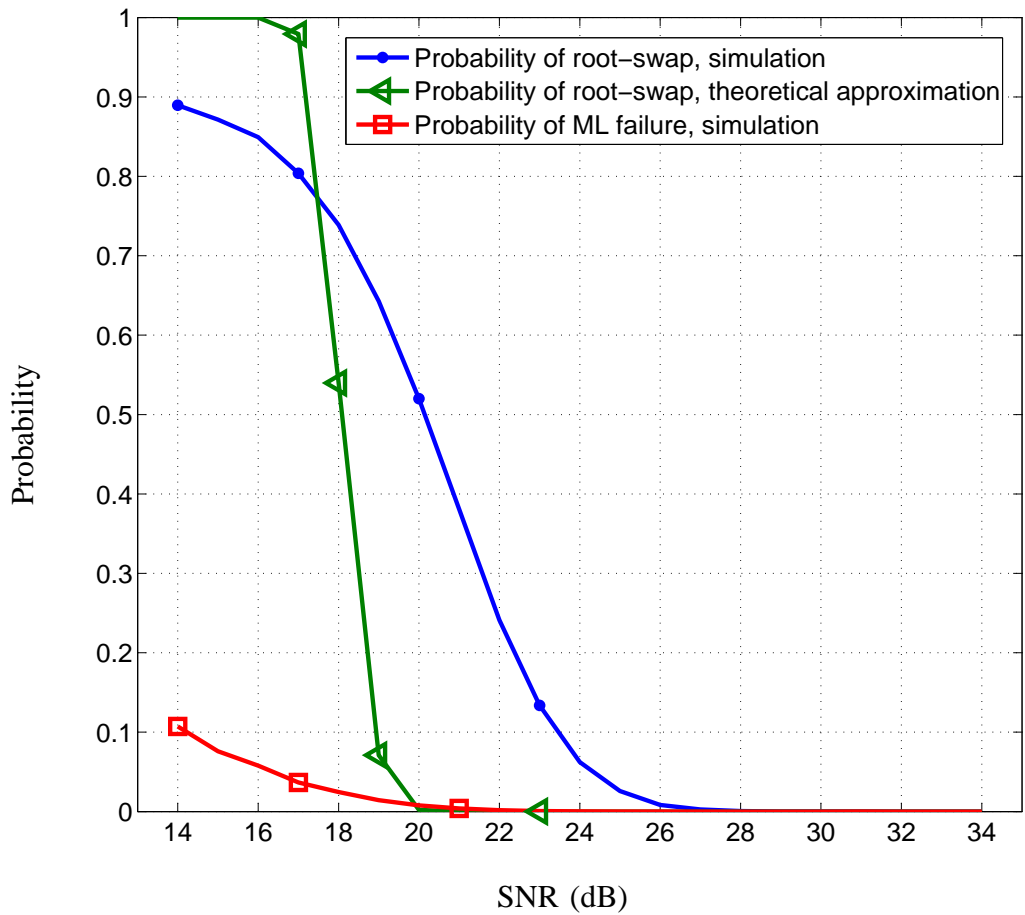}
\vspace{-5mm} \caption{Probability of root-swap and probability of
ML failure versus SNR for correlated sources with $r = 0.9$.
\label{fig:RS_M10_N10_c09}}\vspace{-3mm}
\end{figure}

\begin{figure}[t]
\psfrag{SNR}{SNR (dB)} \psfrag{SL}{Subspace leakage (dB)}
    \includegraphics[width=25cm]{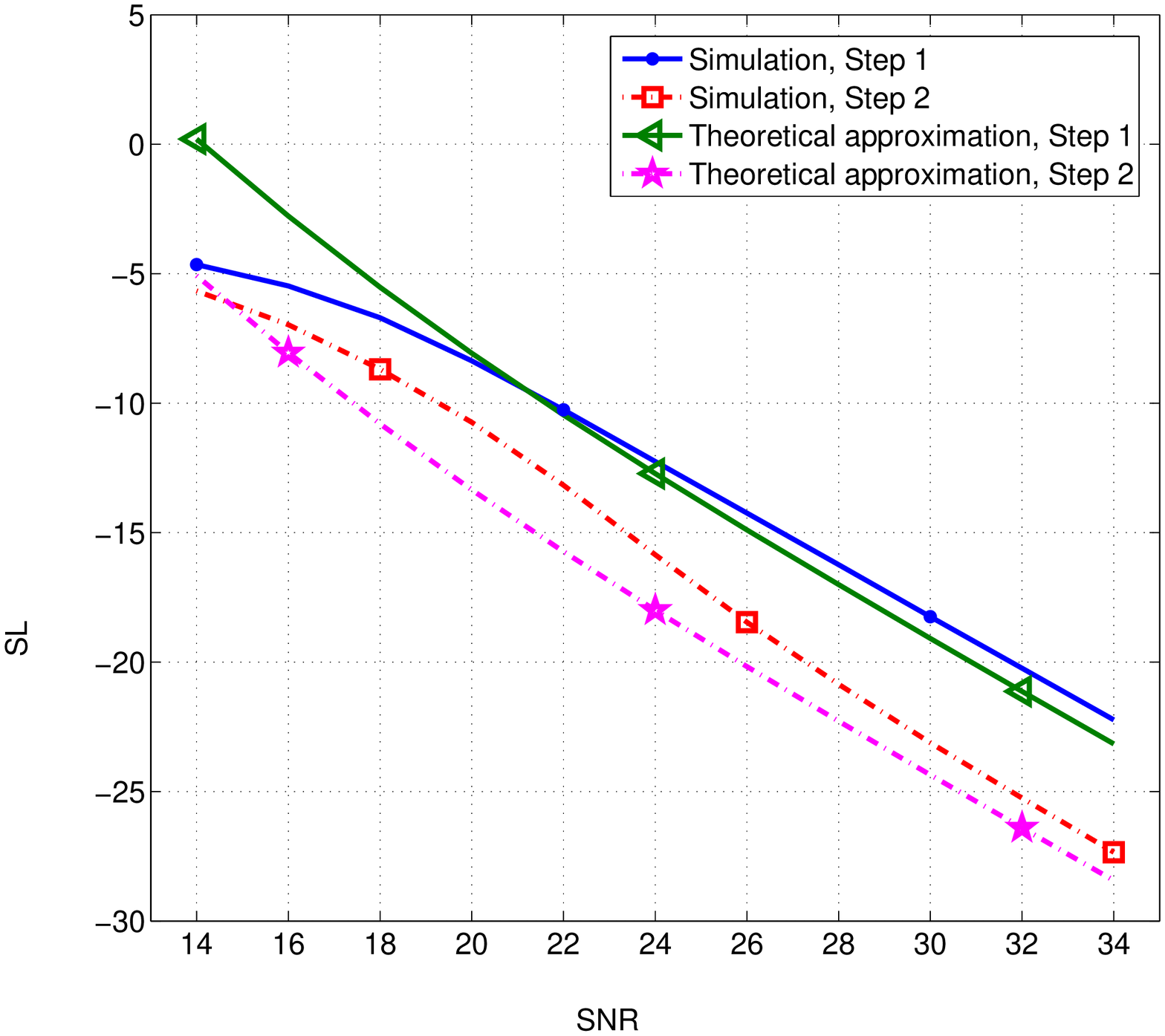}
\vspace{-5mm} \caption{Subspace leakage versus SNR for correlated
sources with $r = 0.9$. The solid and dashed lines represent the
subspace leakage at the first and second steps of the proposed
two-step R-MUSIC algorithm, respectively.
\label{fig:SL_M10_N10_c09}}\vspace{-3mm}
\end{figure}

\begin{figure}[t]
\psfrag{SNR}{SNR (dB)} \psfrag{MSE}{MSE (dB)}
    \includegraphics[width=25cm]{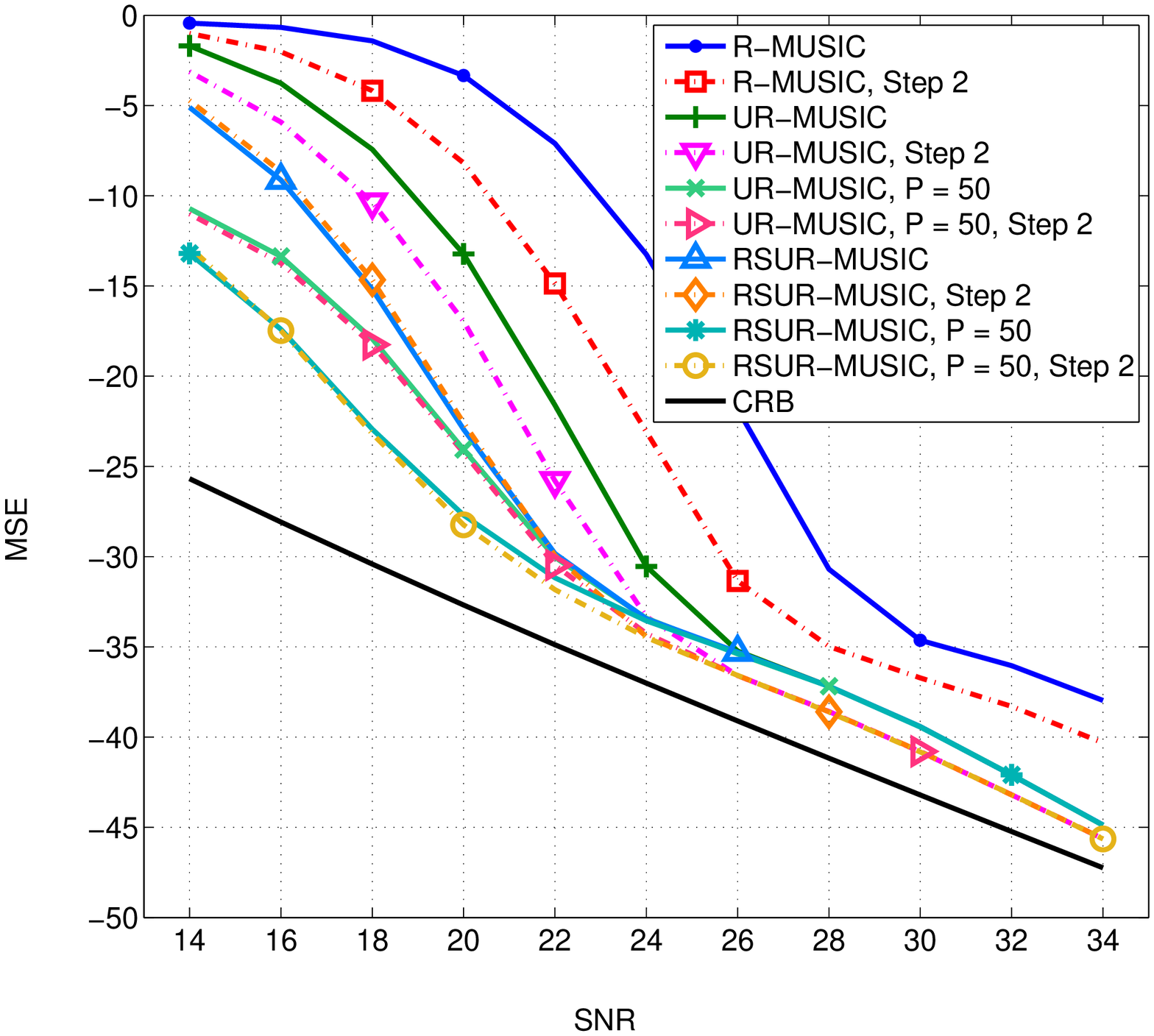}
\vspace{-5mm} \caption{MSE versus SNR for correlated sources with $r
= 0.9$. The solid and dashed lines are based on the first and second
steps of the proposed two-step method, respectively. The methods
used in the two-step algorithm are R-MUSIC, UR-MUSIC, and RSUR-MUSIC
methods. \label{fig:MSE_M10_N10_c09}}\vspace{-3mm}
\end{figure}

\begin{figure}[t]
\psfrag{SNR}{SNR (dB)} \psfrag{PD}{Probability of source resolution}
    \includegraphics[width=25cm]{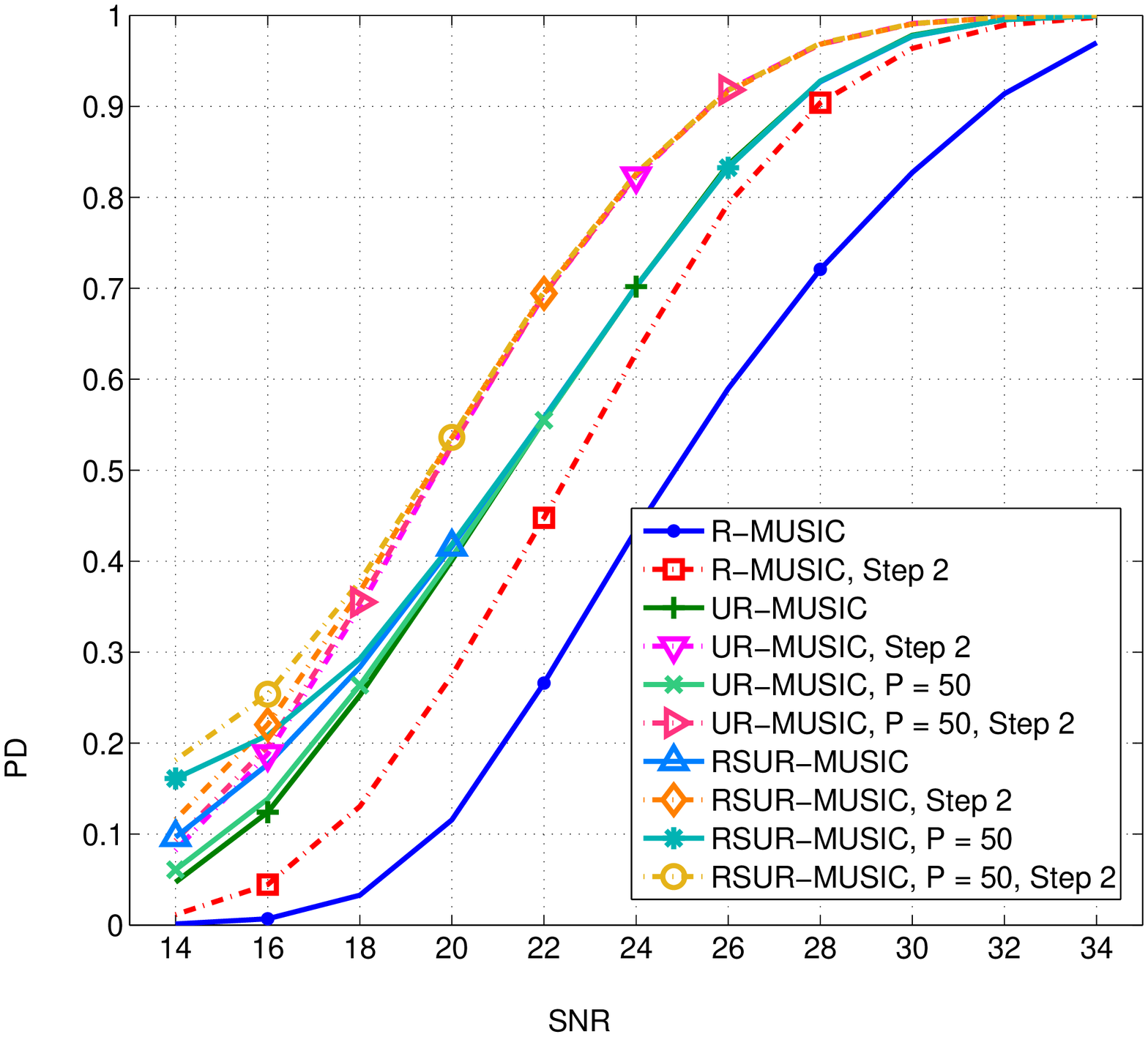}
\vspace{-5mm} \caption{Probability of source resolution versus SNR
for correlated sources  with $r = 0.9$. The solid and dashed lines
are based on the first and second steps of the proposed two-step
method, respectively. The methods used in the two-step algorithm are
R-MUSIC, UR-MUSIC, and RSUR-MUSIC methods.
\label{fig:PD_M10_N10_c09}}\vspace{-3mm}
\end{figure}

\begin{figure}[t]
\psfrag{SNR}{SNR (dB)} \psfrag{CMSE}{CMSE (dB)}
    \includegraphics[width=25cm]{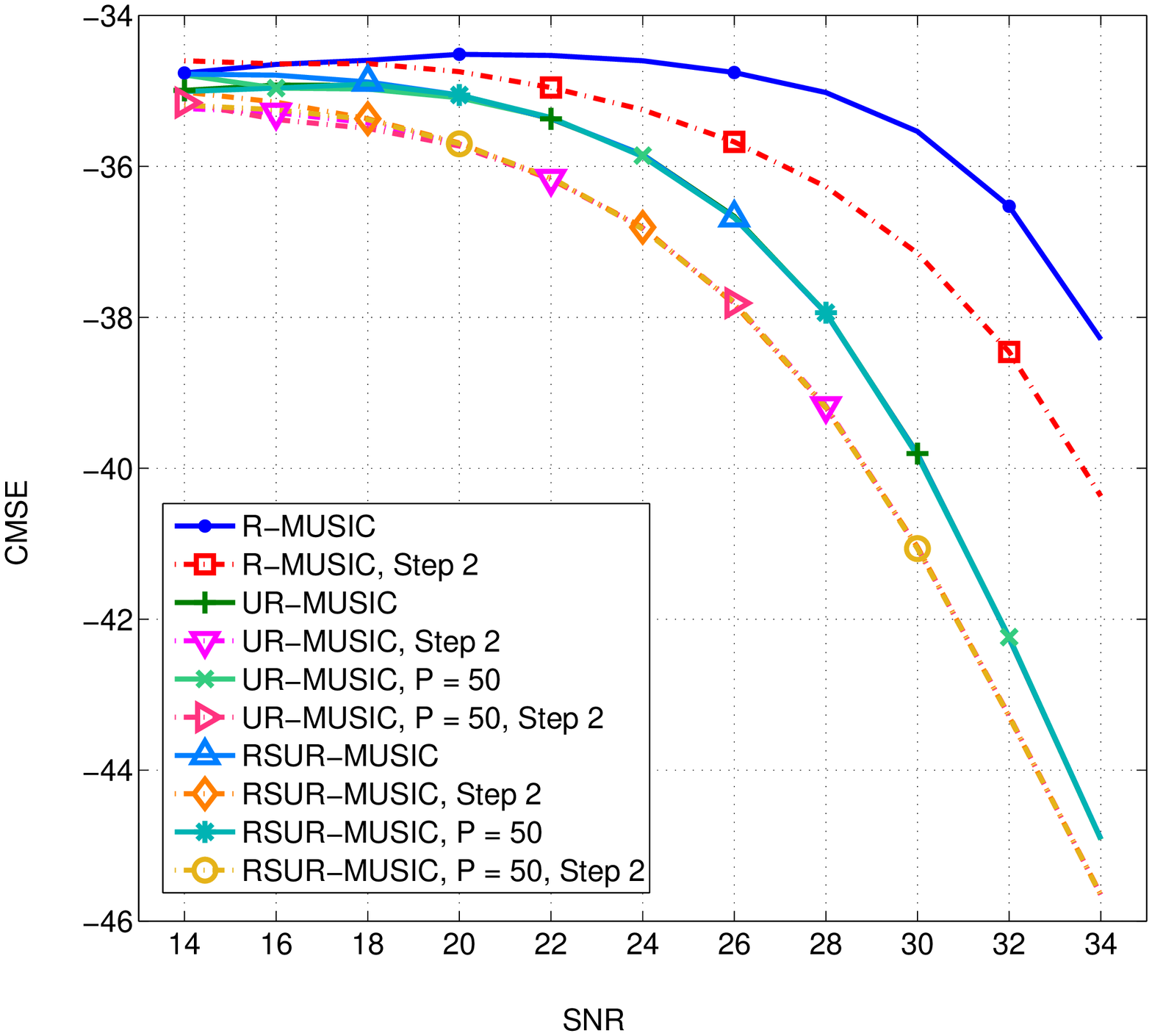}
\vspace{-5mm} \caption{CMSE versus SNR for the correlated sources
with $r = 0.9$. The solid and dashed lines are based on the first
and second steps of the proposed two-step method, respectively. The
methods used in the two-step algorithm are R-MUSIC), UR-MUSIC, and
RSUR-MUSIC methods. \label{fig:CMSE_M10_N10_c09}}\vspace{-3mm}
\end{figure}

The results for the case of the correlated sources with $r = 0.9$
are depicted in Figs.~\ref{fig:RS_M10_N10_c09} to
\ref{fig:CMSE_M10_N10_c09}. Similar observations are made from these
figures as those discussed for the case of the uncorrelated sources.
Compared to the uncorrelated case, the performance breakdown occurs
at a higher SNR value. This makes the importance of the improved
methods more significant, as there is a higher chance that the
actual SNR of a system falls in the breakdown region. As seen from
the figures for the correlated sources, the proposed methods prove
to be helpful in dealing with the performance breakdown problem. The
gain obtained by the improved methods is also more significant
compared to the case of the uncorrelated sources. For instance, the
MSE improvement achieved by the two-step root-MUSIC method for the
uncorrelated sources is about half a dB, while in the case of the
correlated sources, the MSE curve is shifted by more than $2$~dB to
the left. Similarly, more significant performance gains are obtained
for the probability of source resolution and also the CMSE.

\section{Conclusion}\label{sec:conclude}
The performance breakdown of the subspace based DOA estimation
methods in the threshold region where the SNR and/or sample size is
low has been studied in this paper. The subspace leakage as the main
cause of the performance breakdown was formally defined and
theoretically derived. The two-step algorithm has been proposed in
order to reduce the amount of subspace leakage. The introduced
method is based on estimating the DOAs at the first step and
modifying the covariance matrix using the estimated DOAs at the
second step. We have theoretically derived the subspace leakage at
both steps, and have shown that the subspace leakage is reduced at
the second step of the proposed method leading to better
performance. The algorithm can also be extended to the third step by
further modifying the covariance matrix based on the improved
estimates obtained at the second step. We have investigated the
performance of the algorithm for further steps through simulations
(not included in the paper). However, the achieved improvement is
marginal and does not justify the added complexity.
The behavior of the root-MUSIC algorithm in the threshold region has
been also studied, and a phenomenon called root-swap has been
observed to contribute to the performance breakdown. Then, an
improved method has been introduced to remedy this problem by
considering different combinations of the roots and picking up the
one that minimizes the SML function.
The performance improvement achieved by the proposed methods has
also been demonstrated using numerical examples and simulation
results. We also combined the proposed algorithms with the
previously introduced methods in the literature, which resulted in
further improvement in the performance.

\appendices
\section{Probability of Root-Swap Approximation}
\label{sec:appndx_root_swap_prob}
The root-swap is defined as the event when at least one of the
estimated signal roots $\hat{z}_k$ ($1 \leq k \leq K$) has a smaller
magnitude than the magnitude of an estimated noise root $\hat{z}_m$
($K+1 \leq m \leq M-1$), i.e., $\hat{r}_k < \hat{r}_m$. Let us
denote the probability of the event that $\hat{r}_k < \hat{r}_m$ by
$p_{km}$. The complement of this event represents the case when the
$k$-th estimated signal root has not been swapped with the $m$-th
estimated noise root, and its probability is given by $1-p_{km}$.
Let us denote the probability of root-swap by $P(\text{root-swap})$.
The complement of the root-swap event is the event when none of the
estimated signal roots has been swapped with an estimated noise
root, and its probability is given by $1-P(\text{root-swap})$.
Assuming that the individual root-swap events are independent from
each other, we have
\begin{equation}
1-P(\text{root-swap}) = \prod_{k=1}^K \prod_{m=K+1}^{M-1} \left( 1 -
p_{km} \right). \label{eq:root_swap_comp_prob}
\end{equation}

In the sequel, we derive an approximation for $p_{km}$. Noting that
$r_k = 1$ for the true signal roots, we have
\begin{eqnarray}
p_{km} \hspace{-2mm} &=& \hspace{-2mm} P \left( \hat{r}_m >
\hat{r}_k \right) \nonumber \\
&=& \hspace{-2mm} P \left( \Delta r_m - \Delta r_k > 1 - r_m
\right). \label{eq:indiv_root_swap_prob}
\end{eqnarray}

In order to proceed with the computation of $p_{km}$, we consider
the distributions of $\Delta r_m$ and $\Delta r_k$. It is shown in
\cite{Proakis92:root_Music_performance} that $\Delta r_k$ ($1 \leq
k \leq K$) follows the $- \left( \sigma_k / \sqrt{2} \right) \sqrt{
\chi^2\left( 2(M-K)-1\right)}$ distribution where $\chi^2\left( \ell
\right)$ denotes a chi-square distribution with $\ell$ degrees of
freedom and $\sigma_k^2$ is given by
\begin{equation}
\sigma_k^2 = \frac{\sigma_{\text{n}}^2}{N \left( \ba_k^{(1)H}
\bP^{\bot} \ba_k^{(1)} \right) } \sum_{i = 1}^K \frac{\lambda_{M - K
+ i}}{\left( \lambda_{M - K + i} - \sigma_{\text{n}}^2 \right)^2}
\left| \be_i^H \ba_k \right|^2
\end{equation}
where $\bP^{\bot}$ is the true projection matrix of the noise
subspace and $\ba_k^{(1)}$ is given by \eqref{eq:a1_k}.

We next consider the distribution of $\Delta r_m$. In
\cite{Proakis92:root_Music_performance}, the distribution of $\Delta
r_k$ is computed using a second order Taylor expansion of the
estimated root-MUSIC polynomial around the true signal roots (which
are located on the unit circle). The computation of the distribution
of $\Delta r_m$ requires the analysis to be performed around the true
noise roots which are located inside the unit circle. The second
order expansions of $\ba(\hat{z}_k)$ and $\ba^T(\hat{z}_k^{-1})$
around the true signal root $z_k$ are given by
\cite{Proakis92:root_Music_performance}
\begin{eqnarray}
\ba(\hat{z}_k) \hspace{-2mm} & \approx & \hspace{-2mm} \ba_k + j
\ba_k^{(1)} \Delta \omega_k + \ba_k^{(1)} \Delta r_k \nonumber \\
\ba^T(\hat{z}_k^{-1}) \hspace{-2mm} & \approx & \hspace{-2mm}
\ba_k^H - j \ba_k^{(1)H} \Delta \omega_k - \ba_k^{(1)H} \Delta r_k
\label{eq:2_order_taylor_steer_vec_signal}
\end{eqnarray}
where $\Delta \omega_k$ is the difference between the angle of the
$k$-th estimated root and the angle of the corresponding true root.
For the $m$-th noise root, let $\ba_m$  be defined as
\begin{equation}
\ba_m \triangleq \left[1,~e^{-j
\omega_m},\cdots,~e^{-j(M-1)\omega_m} \right]^T
\end{equation}
where $\omega_m$ is the angle of $z_m$. Let also $\ba_m^{(1)}$ be
defined similar to \eqref{eq:a1_k} with $\omega_k$ replaced with
$\omega_m$. Then, the second order expansions of $\ba(\hat{z}_m)$
and $\ba^T(\hat{z}_m^{-1})$ around the true noise root $z_m$ are
given by
\begin{eqnarray}
\hspace{-5mm} \ba(\hat{z}_m) \hspace{-3mm} & \approx & \hspace{-3mm}
\bR_m^{-1} \left( \ba_m + j
\ba_m^{(1)} \Delta \omega_m + \ba_m^{(1)} \left(\frac{\Delta r_m}{r_m}
\right) \right) \nonumber \\
\hspace{-5mm} \ba^T(\hat{z}_m^{-1}) \hspace{-3mm} & \approx &
\hspace{-3mm} \left( \ba_m^H - j \ba_m^{(1)H} \Delta \omega_m -
\ba_m^{(1)H} \left(\frac{\Delta r_m}{r_m}\right) \right)
\hspace{-1mm} \bR_m \label{eq:2_order_taylor_steer_vec_noise}
\end{eqnarray}
where $\bR_m$ is a $M \times M$ diagonal matrix with its diagonal
elements equal to $1,~r_m,\cdots,~r_m^{(M-1)}$. Since the Taylor
expansion for the steering vectors of the roots on the circle and
the expansion for the roots inside the circle, i.e.,
\eqref{eq:2_order_taylor_steer_vec_signal} and
\eqref{eq:2_order_taylor_steer_vec_noise} have similar structures,
it is reasonable to assume that $\Delta r_k$ and $\Delta r_m / r_m$
also have similar distributions. Then, the variance of $\Delta r_m$ is
in the order of the variance of $\Delta r_k$ multiplied by $r_m^2$.
Since $r_m < 1$, the variance of $\Delta r_m$ is smaller than the
variance of $\Delta r_k$. In order to simplify the computation of
$p_{km}$, we ignore the effect of $\Delta r_m$ and approximate
$p_{km}$ by
\begin{equation}
p_{km} \approx P \left( - \Delta r_k > 1 - r_m \right).
\label{eq:indiv_root_swap_prob_approx}
\end{equation}
This is equivalent to using the probability $P \left( \hat{r}_k <
r_m \right)$ as an approximation for $p_{km}$. Since we have the
distribution of $\Delta r_k$, we can compute $p_{km}$ using
\eqref{eq:indiv_root_swap_prob_approx}. When $M - K \gg 1$,
$\Delta r_k$ follows approximately a normal distribution
$\mathcal{N}\left( -\sigma_k \sqrt{M-K-(3/4)}, \sigma_k^2 / 4
\right)$ \cite{Proakis92:root_Music_performance}. Using
\eqref{eq:indiv_root_swap_prob_approx}, the probability $p_{km}$
can be approximated by
\begin{equation}
p_{km} \approx Q\left( \frac{1- r_m - \sigma_k \sqrt{M-K-(3/4)}
}{\sqrt{\sigma_k^2 / 4}} \right).
\label{eq:indiv_root_swap_prob_approx_Gaussian}
\end{equation}
Finally, the approximation of the probability of root-swap
$P(\text{root-swap})$ is found by using the approximation
\eqref{eq:indiv_root_swap_prob_approx_Gaussian}, the expression
\eqref{eq:root_swap_comp_prob}, and the fact that $Q(-x) = 1
- Q(x)$ as
\begin{eqnarray}
&& P(\text{root-swap}) \approx 1 - \prod_{k=1}^K \prod_{m=K+1}^{M-1}
Q\left( \frac{-1 + r_m + \sigma_k \sqrt{M-K-(3/4)}
}{\sqrt{\sigma_k^2 / 4}} \right).
\end{eqnarray}
It completes the derivation.

\section{Subspace Leakage at the First Step}
\label{sec:appndx_SL1}
Let us start with the computation of $\rho_1$. Let $\Delta \bP
\triangleq \widehat{\bP} - \bP$ be the estimation error of the
signal projection matrix. Then, using the properties that $\bP^2 =
\bP$ and $\text{Tr} \left\{ \bP \right\} = K$, the expression
\eqref{eq:sub_leak_def} for the first step of the two-step
root-MUSIC algorithm can be written as
\begin{eqnarray}
\rho_1 \hspace{-2mm} & = & \hspace{-2mm} 1 - \frac{1}{K} \text{Tr}
\left\{ \left( \bP + \Delta \bP \right)
\bP \right\} \nonumber \\
& = & \hspace{-2mm} 1 - \frac{1}{K} \left( K + \text{Tr} \left\{
\Delta \bP
\bP \right\} \right) \nonumber \\
& = & \hspace{-2mm} - \frac{1}{K} \text{Tr} \left\{ \Delta \bP \bP
\right\}.\label{eq:rho1_Tr_DeltaPP}
\end{eqnarray}

It is shown in \cite{Proakis92:root_Music_performance} that the
series expansion of $\widehat{\bP}$ based on $\Delta \bR$ is given
by
\begin{equation}
\widehat{\bP} = \bP + \delta \bP + \cdots + \delta^n \bP + \cdots
\label{eq:Phat_series_exp}
\end{equation}
where
\begin{equation}
\delta \bP = \bP^{\bot} \Delta \bR \bV^{\dag} + \bV^{\dag} \Delta
\bR \bP^{\bot} \label{eq:Phat_series_exp_deltaP}
\end{equation}
and the rest of the terms are related by the following recurrence
\begin{eqnarray}
\delta^n \bP \hspace{-2mm} & = & \hspace{-2mm} -\bP^{\bot} \left(
\delta^{n-1} \bP \right) \Delta \bR \bV^{\dag} + \bP^{\bot} \Delta
\bR \left( \delta^{n-1} \bP \right)
\bV^{\dag} \nonumber \\
&& \hspace{-2mm} - \bV^{\dag} \Delta \bR \left( \delta^{n-1} \bP
\right) \bP^{\bot} + \bV^{\dag} \left( \delta^{n-1} \bP \right)
\Delta \bR \bP^{\bot} \nonumber \\
&& \hspace{-2mm} - \sum_{i = 1}^{n - 1} \bP \left( \delta^i \bP
\right) \left( \delta^{n-i} \bP \right) \bP \nonumber \\
&& \hspace{-2mm} + \sum_{i = 1}^{n - 1} \bP^{\bot} \left( \delta^i
\bP \right) \left( \delta^{n-i} \bP \right) \bP^{\bot}.
\label{eq:Phat_series_exp_deltanP}
\end{eqnarray}

The following lemma will be further used.

\theoremstyle{plain}
\newtheorem{lemma_PVdag}{Lemma}[]
\begin{lemma_PVdag}\label{lem:PVdag}
The columns of $\bV^{\dag}$ belong to the signal subspace, i.e.,
$\bP \bV^{\dag} = \bV^{\dag}$.
\end{lemma_PVdag}
\begin{proof}
The proof follows by multiplying $\bP$ by $\bV^{\dag}$ and then
substituting $\bP$ with $\bE \bE^H$ and $\bV^{\dag}$ with
\eqref{eq:pseudoinverse_V}.
\end{proof}

In a similar way to Lemma~\ref{lem:PVdag}, it can also be shown that
\begin{equation}
\bV \bV^{\dag} = \bV^{\dag} \bV = \bP. \label{eq:VVdag_VdagV_P}
\end{equation}

Using \eqref{eq:rho1_Tr_DeltaPP}, the series expansion of
$\widehat{\bP}$ in \eqref{eq:Phat_series_exp}, expressions
\eqref{eq:Phat_series_exp_deltaP} and
\eqref{eq:Phat_series_exp_deltanP} up to the $\delta^2 \bP$ term,
and the facts that $\bP \bP^{\bot} = \bP^{\bot} \bP =
\boldsymbol{0}$ and $\bP \bP = \bP$, we can write $\rho_1$ as
\begin{equation}
\rho_1 = - \frac{1}{K} \text{Tr} \left\{ -\bP \left( \delta \bP
\right) \left( \delta \bP \right)
\right\}.\label{eq:rho1_Tr_PDeltaPDeltaP}
\end{equation}
Then, $\rho_1$ is computed by substituting
\eqref{eq:Phat_series_exp_deltaP} in
\eqref{eq:rho1_Tr_PDeltaPDeltaP}, using $\bP^{\bot} \bP^{\bot} =
\bP^{\bot}$, and Lemma~\ref{lem:PVdag} as
\begin{eqnarray}
\rho_1 \hspace{-2mm} & = & \hspace{-2mm} \frac{1}{K} \text{Tr}
\Big\{ \bP \left( \bP^{\bot} \Delta \bR \bV^{\dag} + \bV^{\dag}
\Delta \bR \bP^{\bot} \right) \left( \bP^{\bot} \Delta \bR
\bV^{\dag} + \bV^{\dag} \Delta
\bR \bP^{\bot} \right) \Big\} \nonumber \\
& = & \hspace{-2mm} \frac{1}{K} \text{Tr} \left\{ \bP \bV^{\dag}
\Delta \bR \bP^{\bot} \bP^{\bot} \Delta \bR \bV^{\dag} \right\} \nonumber \\
& = & \hspace{-2mm} \frac{1}{K} \text{Tr} \left\{ \bV^{\dag} \Delta
\bR \bP^{\bot} \Delta \bR \bV^{\dag}
\right\}.\label{eq:rho1_Tr_VdagDeltaRPbot}
\end{eqnarray}

Computation of the expected value of the subspace leakage requires
considering the statistical properties of $\Delta \bR$. We use the
following two properties in our derivations
\cite{Proakis92:root_Music_performance}.

\theoremstyle{plain}
\newtheorem{lemma_EDeltaR}[lemma_PVdag]{Lemma}
\begin{lemma_EDeltaR}\label{lem:EDeltaR}
For all matrices $\bA_1,~\bA_2 \in \mathbb{C}^{M \times M}$, we have
\begin{equation}
E\left\{ \Delta \bR \bA_1 \Delta \bR \right\} = \frac{1}{N}
\text{Tr} \left\{ \bR \bA_1 \right\} \bR \label{eq:DeltaR_lemma_P1}
\end{equation}
and
\begin{equation}
E\big\{ \text{Tr} \left\{ \Delta \bR \bA_1 \right\} \text{Tr}
\left\{ \Delta \bR \bA_2 \right\} \big\} = \frac{1}{N} \text{Tr}
\left\{ \bR \bA_1 \bR \bA_2 \right\}. \label{eq:DeltaR_lemma_P2}
\end{equation}
\end{lemma_EDeltaR}

Using \eqref{eq:rho1_Tr_VdagDeltaRPbot} and
\eqref{eq:DeltaR_lemma_P1}, the expected value of $\rho_1$ can be
computed as
\begin{eqnarray}
E\left\{ \rho_1 \right\} \hspace{-2mm} & = & \hspace{-2mm}
\frac{1}{K} \text{Tr} \left\{ \bV^{\dag} E\left\{ \Delta \bR
\bP^{\bot} \Delta \bR \right\} \bV^{\dag} \right\} \nonumber \\
& = & \hspace{-2mm} \frac{1}{K} \text{Tr} \left\{ \bV^{\dag}
\frac{1}{N} \text{Tr} \left\{ \bR \bP^{\bot} \right\} \bR \bV^{\dag}
\right\} \nonumber \\
& = & \hspace{-2mm} \frac{1}{NK} \text{Tr} \left\{ \bP^{\bot} \bR
\right\} \text{Tr} \left\{ \bV^{\dag} \bV^{\dag} \bR \right\}.
\label{eq:E_rho1_Tr_PbotR_Tr_VdagVdagR}
\end{eqnarray}
Since the range space of the matrix $\bA$ is the same as the signal
subspace, we have $\bP^{\bot} \bA = \boldsymbol{0}$. As a result,
$\text{Tr} \left\{ \bP^{\bot} \bR \right\}$ can be simplified as
\begin{eqnarray}
\text{Tr} \left\{ \bP^{\bot} \bR \right\} \hspace{-2mm} & = &
\hspace{-2mm} \text{Tr} \left\{ \bP^{\bot} \left( \bA\bS\bA^H
+ \sigma_{\text{n}}^2\bI_M \right) \right\} \nonumber \\
& = & \hspace{-2mm} \text{Tr} \left\{ \sigma_{\text{n}}^2 \bP^{\bot}
\right\} = \sigma_{\text{n}}^2 \text{Tr} \left\{ \bI_M -
\bP \right\} \nonumber \\
& = & \hspace{-2mm} \sigma_{\text{n}}^2 \left( M - K \right).
\label{eq:Tr_PbotR}
\end{eqnarray}
Furthermore, using \eqref{eq:pseudoinverse_V} and the fact that the
eigenvectors of $\bR$ are orthonormal, the product $\bV^{\dag}
\bV^{\dag} \bR$ can be written as
\begin{equation}
\bV^{\dag} \bV^{\dag} \bR = \sum_{k = 1}^K \frac{\lambda_{M - K +
k}}{\left( \lambda_{M - K + k} - \sigma_{\text{n}}^2 \right)^2 }
\be_k \be_k^H
\end{equation}
which results in
\begin{equation}
\text{Tr} \left\{ \bV^{\dag} \bV^{\dag} \bR \right\} = \sum_{k =
1}^K \frac{\lambda_{M - K + k}}{\left( \lambda_{M - K + k} -
\sigma_{\text{n}}^2 \right)^2 }. \label{eq:Tr_VdagVdagR}
\end{equation}
Finally, $E\left\{ \rho_1 \right\}$ is obtained by substituting
\eqref{eq:Tr_PbotR} and \eqref{eq:Tr_VdagVdagR} in
\eqref{eq:E_rho1_Tr_PbotR_Tr_VdagVdagR} as
\begin{equation}
E\left\{ \rho_1 \right\} = \frac{\sigma_{\text{n}}^2 \left( M - K
\right)}{NK} \sum_{k = 1}^K \frac{\lambda_{M - K + k}}{\left(
\lambda_{M - K + k} - \sigma_{\text{n}}^2 \right)^2 }.
\end{equation}

\section{Subspace Leakage at the Second Step}\label{sec:appndx_SL2}
The subspace leakage at the second step of the two-step root-MUSIC
algorithm can be obtained through the same steps taken for the
computation of $\rho_1$. Referring to
\eqref{eq:rho1_Tr_VdagDeltaRPbot}, the subspace leakage $\rho_2$ is
given by
\begin{equation}
\rho_2 = \frac{1}{K} \text{Tr} \left\{ \bV^{\dag} \Delta \bR^{(2)}
\bP^{\bot} \Delta \bR^{(2)} \bV^{\dag} \right\}
\label{eq:rho2_Tr_VdagDeltaRPbot}
\end{equation}
where $\Delta \bR^{(2)} \triangleq \widehat{\bR}^{(2)} - \bR$ is the
estimation error of the covariance matrix at the second step of the
algorithm. Using \eqref{eq:step2_Rhat2}, the estimation error
$\Delta \bR^{(2)}$ is given by
\begin{equation}
\Delta \bR^{(2)} = \Delta \bR - \gamma \left( \bT+ \bT^H \right).
\label{eq:DeltaR_2}
\end{equation}
Recalling \eqref{eq:T_PARPAbot}, we have $\bT = \widehat{\bP}_A
\widehat{\bR} \widehat{\bP}_A^{\bot}$.

Consider the first order Taylor series expansion of
$\widehat{\bP}_A$ around the true DOAs given by
\begin{equation}
\widehat{\bP}_A \approx \bP_A + d\bP \label{eq:Taylor1_PAhat}
\end{equation}
where $\bP_A \triangleq \bA \left( \bA^H \bA \right)^{-1} \bA^H$ is
equal to the true signal projection matrix\footnote{Note that
although $\bP_A$ is equal to $\bP$, the estimates $\widehat{\bP}_A$
and $\widehat{\bP}$ are obtained in different ways and are not
essentially equal to each other.}, i.e., $\bP_A = \bP$, and $d\bP$
is given by
\begin{equation}
d\bP = \sum_{k = 1}^K \frac{\partial \bP_A}{\partial \omega_k}
\Delta \omega_k. \label{eq:dP_partialP_Deltaomega}
\end{equation}
Here $\Delta \omega_k \triangleq \hat{\omega}_k - \omega_k$ is the
estimation error of $\omega_k$ with $\hat{\omega}_k \triangleq 2 \pi
(d/\lambda) \sin(\hat{\theta}_k)$.

Note that for any square and invertible matrix $\bB$, the partial
derivative of $\bB^{-1}$ with respect to the variable $\omega$ is
given by \cite{Peterson12:matrix_cookbook}
\begin{equation}
\frac{\partial \bB^{-1}}{\partial \omega} = - \bB^{-1}
\frac{\partial \bB}{\partial \omega}
\bB^{-1}.\label{eq:lemma_mat_inv_diff}
\end{equation}
Using \eqref{eq:lemma_mat_inv_diff}, the partial derivative
$\partial \bP_A /
\partial \omega_k$ can be computed as
\begin{eqnarray}
\frac{\partial \bP_A}{\partial \omega_k} \hspace{-2mm} & = &
\hspace{-2mm} \frac{\partial \bA}{\partial \omega_k} \left( \bA^H
\bA \right)^{-1} \bA^H + \bA \frac{\partial \left( \bA^H \bA
\right)^{-1}}{\partial \omega_k} \bA^H  + \bA \left( \bA^H \bA
\right)^{-1} \left( \frac{\partial \bA}{\partial \omega_k}
\right)^H \nonumber \\
& = & \hspace{-2mm} \frac{\partial \bA}{\partial \omega_k} \left(
\bA^H \bA \right)^{-1} \bA^H - \bA \left( \bA^H \bA \right)^{-1}
\left( \left( \frac{\partial \bA}{\partial \omega_k} \right)^H \bA +
\bA^H \frac{\partial \bA}{\partial \omega_k} \right) \left( \bA^H
\bA
\right)^{-1} \bA^H \nonumber \\
&& \hspace{-2mm} + \bA \left( \bA^H \bA \right)^{-1} \left(
\frac{\partial \bA}{\partial \omega_k} \right)^H.
\end{eqnarray}
Then, using \eqref{eq:P_orth+P=I} and $\bP = \bA \left( \bA^H \bA
\right)^{-1} \bA^H$, the partial derivative $\partial \bP_A /
\partial \omega_k$ is given by
\begin{eqnarray}
\frac{\partial \bP_A}{\partial \omega_k} \hspace{-2mm} & = &
\hspace{-2mm} \bP^{\bot} \frac{\partial \bA}{\partial \omega_k}
\left( \bA^H \bA \right)^{-1} \bA^H + \bA \left( \bA^H \bA
\right)^{-1} \left( \frac{\partial \bA}{\partial \omega_k} \right)^H
\bP^{\bot}. \label{eq:partial_PA_partial_omegak}
\end{eqnarray}

The estimation error of $\omega_k$, i.e., $\Delta \omega_k$ in
\eqref{eq:dP_partialP_Deltaomega}, can be written based on $\Delta
\bR$ as \cite{Proakis92:root_Music_performance}
\begin{equation}
\Delta \omega_k = \frac{\ba_k^{(1)H} \bP^{\bot} \Delta \bR
\bV^{\dag} \ba_k - \ba_k^H \bV^{\dag} \Delta \bR \bP^{\bot}
\ba_k^{(1)}}{2 j \left( \ba_k^{(1)H} \bP^{\bot} \ba_k^{(1)}
\right)}. \label{eq:Delta_omega_k}
\end{equation}

The first order Taylor series expansion of $\widehat{\bP}_A^{\bot}$
is obtained using \eqref{eq:PAhatbot_I_PAhat} and
\eqref{eq:Taylor1_PAhat} as
\begin{equation}
\widehat{\bP}_A^{\bot} \approx \bP_A^{\bot} - d\bP
\label{eq:Taylor1_PAhatbot}
\end{equation}
where $\bP_A^{\bot} \triangleq \bI_M - \bP_A$.

The matrix $\bT$ can be then computed using expressions
\eqref{eq:T_PARPAbot}, \eqref{eq:Taylor1_PAhat}, and
\eqref{eq:Taylor1_PAhatbot} with keeping only the first order terms
and noting that $\bP_A = \bP$, $\bP_A^{\bot} = \bP^{\bot}$, and $\bP
\bR \bP^{\bot} = \boldsymbol{0}$ as
\begin{eqnarray}
\bT \hspace{-2mm} & = & \hspace{-2mm} \left( \bP_A + d\bP \right)
\left( \bR + \Delta R \right) \left( \bP_A^{\bot} - d\bP \right)
\nonumber \\
& \approx & \hspace{-2mm} - \bP \bR d\bP + \bP \Delta R \bP^{\bot} +
d\bP \bR \bP^{\bot}. \label{eq:T_approx_first_order}
\end{eqnarray}

We can now compute $\rho_2$ using expressions
\eqref{eq:rho2_Tr_VdagDeltaRPbot}, \eqref{eq:DeltaR_2}, and
\eqref{eq:T_approx_first_order} as
\begin{eqnarray}
\rho_2 \hspace{-2mm} & = & \hspace{-2mm} \frac{1}{K} \text{Tr}
\Big\{ \bV^{\dag} \left( \Delta \bR - \gamma \left( \bT+ \bT^H
\right) \right) \bP^{\bot} \left( \Delta \bR - \gamma \left( \bT+
\bT^H \right) \right) \bV^{\dag} \Big\} \nonumber \\
& = & \hspace{-2mm} \frac{1}{K} \text{Tr} \Big\{ \bV^{\dag} \big(
\Delta \bR - \gamma \big( - \bP \bR d\bP + \bP \Delta R \bP^{\bot} +
d\bP \bR \bP^{\bot} - d\bP \bR \bP + \bP^{\bot} \Delta R \bP \nonumber \\
&& \hspace{29mm}  + \bP^{\bot} \bR d\bP \big) \big) \bP^{\bot} \big(
\Delta \bR - \gamma \big( - \bP \bR d\bP + \bP \Delta R \bP^{\bot} + d\bP \bR \bP^{\bot} \nonumber \\
&& \hspace{29mm} - d\bP \bR \bP + \bP^{\bot} \Delta R \bP +
\bP^{\bot} \bR d\bP \big) \big) \bV^{\dag} \Big\}.
\end{eqnarray}
Then, using expressions \eqref{eq:dP_partialP_Deltaomega},
\eqref{eq:partial_PA_partial_omegak}, and the fact that $\bP
\bP^{\bot} = \bP^{\bot} \bP = \bV^{\dag} \bP^{\bot} = \bP^{\bot}
\bV^{\dag} = \boldsymbol{0}$ to eliminate the terms that equal zero,
$\rho_2$ is computed as
\begin{eqnarray}
\rho_2 \hspace{-2mm} & = & \hspace{-2mm} \frac{1}{K} \text{Tr}
\Big\{ \bV^{\dag} \big( \Delta \bR - \gamma \big( - \bP \bR d\bP +
\bP \Delta R \bP^{\bot} + d\bP \bR \bP^{\bot} \big) \big) \nonumber \\
&& \hspace{14mm} \times \bP^{\bot} \big( \Delta \bR - \gamma \big( -
d\bP \bR \bP + \bP^{\bot} \Delta R \bP + \bP^{\bot} \bR d\bP \big)
\big) \bV^{\dag} \Big\}. \label{eq:rho2_Taylor_Expansion}
\end{eqnarray}
Expanding the terms in \eqref{eq:rho2_Taylor_Expansion} and using
the fact that $\bP \bV^{\dag} = \bV^{\dag} \bP = \bV^{\dag}$ results
in the following expression for $\rho_2$
\begin{eqnarray}
\rho_2 \hspace{-3mm} & = & \hspace{-3mm} \frac{1}{K} \text{Tr}
\Big\{ \bV^{\dag} \Delta \bR \bP^{\bot} \Delta \bR \bV^{\dag} -
\hspace{-0.5mm} \gamma \big( \hspace{-0.8mm} - \hspace{-0.8mm}
\bV^{\dag} \Delta \bR \bP^{\bot} d\bP \bR \bV^{\dag} + \bV^{\dag}
\Delta \bR \bP^{\bot} \Delta \bR \bV^{\dag} \hspace{-0.5mm}
\hspace{-0.5mm} \nonumber \\
&& \hspace{4mm}  + \bV^{\dag} \Delta \bR \bP^{\bot} \bR d\bP
\bV^{\dag} - \bV^{\dag} \bR d\bP \bP^{\bot} \Delta \bR \bV^{\dag} +
\bV^{\dag}
\Delta \bR \bP^{\bot} \Delta \bR \bV^{\dag} \nonumber \\
&& \hspace{4mm} + \bV^{\dag} d\bP \bR \bP^{\bot} \Delta \bR
\bV^{\dag} \big) + \gamma^2 \big( \bV^{\dag} \bR d\bP \bP^{\bot}
d\bP \bR \bV^{\dag} - \bV^{\dag} \bR d\bP \bP^{\bot} \Delta \bR
\bV^{\dag}\nonumber \\
&& \hspace{4mm} - \bV^{\dag} \bR d\bP \bP^{\bot} \bR d\bP \bV^{\dag}
- \bV^{\dag} \Delta \bR \bP^{\bot} d\bP \bR \bV^{\dag} + \bV^{\dag}
\Delta \bR \bP^{\bot} \Delta \bR
\bV^{\dag} \nonumber \\
&& \hspace{4mm} + \bV^{\dag} \Delta \bR \bP^{\bot} \bR d\bP
\bV^{\dag} - \bV^{\dag} d\bP \bR \bP^{\bot} d\bP \bR \bV^{\dag} +
\bV^{\dag} d\bP \bR \bP^{\bot} \Delta \bR
\bV^{\dag} \nonumber \\
&& \hspace{4mm} + \bV^{\dag} d\bP \bR \bP^{\bot} \bR d\bP \bV^{\dag}
\big)
\Big\}.\nonumber \\
\label{eq:rho2_big_expansion}
\end{eqnarray}
By reordering the terms in \eqref{eq:rho2_big_expansion}, the
subspace leakage $\rho_2$ can be further rewritten as
\begin{eqnarray}
\rho_2 \hspace{-3mm} & = & \hspace{-3mm} \frac{1}{K} \text{Tr}
\Big\{ \left( 1 -2\gamma + \gamma^2 \right) \bV^{\dag} \Delta \bR
\bP^{\bot} \Delta \bR \bV^{\dag} + \left( \gamma^2 - \gamma \right)
\big( - \bV^{\dag} \Delta \bR \bP^{\bot} d\bP \bR
\bV^{\dag} \nonumber \\
&& \hspace{4mm} + \bV^{\dag} \Delta \bR \bP^{\bot} \bR d\bP
\bV^{\dag} - \bV^{\dag} \bR d\bP \bP^{\bot} \Delta \bR \bV^{\dag} +
\bV^{\dag} d\bP \bR \bP^{\bot} \Delta \bR \bV^{\dag} \big)
\nonumber \\
&& \hspace{4mm} + \gamma^2 \big( \bV^{\dag} \bR d\bP \bP^{\bot} d\bP
\bR \bV^{\dag} - \bV^{\dag} \bR d\bP \bP^{\bot} \bR d\bP \bV^{\dag}
- \bV^{\dag} d\bP \bR \bP^{\bot} d\bP \bR \bV^{\dag} \nonumber \\
&& \hspace{4mm} + \bV^{\dag} d\bP \bR \bP^{\bot} \bR d\bP \bV^{\dag}
\big) \Big\}.\label{eq:rho2_reordered_expansion}
\end{eqnarray}

The terms multiplied by $\left( \gamma^2 - \gamma \right)$ in
\eqref{eq:rho2_reordered_expansion} can be simplified using
expressions \eqref{eq:V_R_sigma2I}, \eqref{eq:VVdag_VdagV_P}, and
the fact that $\bP^{\bot} \bV = \boldsymbol{0}$ as
\begin{eqnarray}
&& - \bV^{\dag} \Delta \bR \bP^{\bot} d\bP \left( \bV +
\sigma_{\text{n}}^2\bI_M \right) \bV^{\dag} + \bV^{\dag} \Delta \bR
\bP^{\bot} \left( \bV + \sigma_{\text{n}}^2\bI_M \right) d\bP
\bV^{\dag} \nonumber \\
&& - \bV^{\dag} \left( \bV + \sigma_{\text{n}}^2\bI_M \right) d\bP
\bP^{\bot} \Delta \bR \bV^{\dag} + \bV^{\dag} d\bP \left( \bV +
\sigma_{\text{n}}^2\bI_M \right) \bP^{\bot} \Delta \bR \bV^{\dag}
\nonumber \\
&& = - \bV^{\dag} \Delta \bR \bP^{\bot} d\bP \bP - \bP d\bP
\bP^{\bot} \Delta \bR \bV^{\dag}.\label{eq:rho2_gamma2_gamma_terms}
\end{eqnarray}
In a similar way, the terms multiplied by $\gamma^2$ in
\eqref{eq:rho2_reordered_expansion} can be simplified to
\begin{eqnarray}
&& \bV^{\dag} \bR d\bP \bP^{\bot} d\bP \left( \bV +
\sigma_{\text{n}}^2\bI_M \right) \bV^{\dag} - \bV^{\dag} \bR d\bP
\bP^{\bot} \left( \bV + \sigma_{\text{n}}^2\bI_M \right) d\bP \bV^{\dag} \nonumber \\
&& - \bV^{\dag} d\bP \bR \bP^{\bot} d\bP \left( \bV +
\sigma_{\text{n}}^2\bI_M \right) \bV^{\dag} + \bV^{\dag} d\bP \bR
\bP^{\bot} \left( \bV + \sigma_{\text{n}}^2\bI_M \right) d\bP
\bV^{\dag}
\nonumber \\
&& = \bV^{\dag} \bR d\bP \bP^{\bot} d\bP \bP - \bV^{\dag} d\bP \bR
\bP^{\bot} d\bP \bP \nonumber \\
&& = \bV^{\dag} \left( \bV + \sigma_{\text{n}}^2\bI_M \right) d\bP
\bP^{\bot} d\bP \bP - \bV^{\dag} d\bP \left( \bV +
\sigma_{\text{n}}^2\bI_M \right) \bP^{\bot} d\bP \bP \nonumber \\
&& = \bP d\bP \bP^{\bot} d\bP \bP
\end{eqnarray}
which using the fact that $\bP^{\bot} d\bP \bP^{\bot} =
\boldsymbol{0}$ (see \eqref{eq:dP_partialP_Deltaomega} and
\eqref{eq:partial_PA_partial_omegak}) can be further simplified to
\begin{eqnarray}
\bP d\bP \bP^{\bot} d\bP \bP \hspace{-2mm} & = & \hspace{-2mm}
\left( \bI_M -\bP^{\bot} \right) d\bP \bP^{\bot} d\bP \left(
\bI_M -\bP^{\bot} \right) \nonumber \\
& = & \hspace{-2mm} d\bP \bP^{\bot} d\bP.
\label{eq:rho2_gamma2_terms}
\end{eqnarray}

Finally, using expressions \eqref{eq:rho1_Tr_VdagDeltaRPbot},
\eqref{eq:rho2_reordered_expansion},
\eqref{eq:rho2_gamma2_gamma_terms}, \eqref{eq:rho2_gamma2_terms},
and Lemma~\ref{lem:PVdag}, the subspace leakage $\rho_2$ is computed
as
\begin{eqnarray}
\rho_2 \hspace{-2.5mm} & = & \hspace{-2.5mm} \left( 1 - 2\gamma +
\gamma^2 \right) \rho_1 + \frac{ 2 \left( \gamma - \gamma^2 \right)
}{K} Re \left\{ \text{Tr} \left\{ \bV^{\dag} \Delta \bR \bP^{\bot}
d\bP \right\} \right\} + \frac{\gamma^2}{K} \text{Tr} \left\{ d\bP
\bP^{\bot} d\bP \right\}. \label{eq:rho2_final_simplified}
\end{eqnarray}

Computation of the expected value of $\rho_2$ involves finding the
expected value of the two trace functions in
\eqref{eq:rho2_final_simplified}. Using expressions
\eqref{eq:dP_partialP_Deltaomega} and
\eqref{eq:partial_PA_partial_omegak}, the expected value of the
first trace function in \eqref{eq:rho2_final_simplified} is given by
\begin{eqnarray}
E \left\{ \text{Tr} \left\{ \bV^{\dag} \Delta \bR \bP^{\bot} d\bP
\right\} \right\}  = E \left\{ \text{Tr} \left\{ \Delta \bR \sum_{k
= 1}^K \bP^{\bot} \frac{\partial \bA}{\partial \omega_k} \left(
\bA^H \bA \right)^{-1} \bA^H \Delta \omega_k \bV^{\dag} \right\}
\right\}. \label{eq:trace1_0_rho2}
\end{eqnarray}
Then, by substituting \eqref{eq:Delta_omega_k} in
\eqref{eq:trace1_0_rho2}, we have
\begin{eqnarray}
&& \hspace{-16mm} E \left\{ \text{Tr} \left\{ \bV^{\dag} \Delta \bR
\bP^{\bot} d\bP \right\} \right\} = E \Bigg\{ \text{Tr} \Bigg\{
\sum_{k = 1}^K \Delta \bR \bP^{\bot} \frac{\partial \bA}{\partial
\omega_k} \left( \bA^H
\bA \right)^{-1} \bA^H \bV^{\dag} \nonumber \\
&& \hspace{15mm} \times \frac{1}{2 j \left( \ba_k^{(1)H} \bP^{\bot}
\ba_k^{(1)} \right)} \Big( \ba_k^{(1)H} \bP^{\bot} \Delta \bR
\bV^{\dag} \ba_k - \ba_k^H \bV^{\dag} \Delta \bR \bP^{\bot}
\ba_k^{(1)} \Big) \Bigg\} \Bigg\}. \label{eq:trace1_1_rho2}
\end{eqnarray}
The order of the summation and trace operator in
\eqref{eq:trace1_1_rho2} can be swaped. Moreover, the last two terms
can be written using the trace operator as
\begin{eqnarray}
&& \hspace{-15mm} E \left\{ \text{Tr} \left\{ \bV^{\dag} \Delta \bR
\bP^{\bot} d\bP \right\} \right\} = E \Bigg\{ \sum_{k = 1}^K
\frac{1}{2 j \left(
\ba_k^{(1)H} \bP^{\bot} \ba_k^{(1)} \right)} \nonumber \\
&& \hspace{30mm} \times \text{Tr} \Bigg\{ \Delta \bR \bP^{\bot}
\frac{\partial \bA}{\partial \omega_k} \left( \bA^H \bA \right)^{-1}
\bA^H \bV^{\dag} \Bigg\} \Big( \text{Tr} \left\{ \Delta \bR
\bV^{\dag} \ba_k \ba_k^{(1)H} \bP^{\bot} \right\} \nonumber \\
&& \hspace{30mm} - \text{Tr} \left\{ \Delta \bR \bP^{\bot}
\ba_k^{(1)} \ba_k^H \bV^{\dag} \right\} \Big) \Bigg\}.
\label{eq:trace1_2_rho2}
\end{eqnarray}
The expression in \eqref{eq:trace1_2_rho2} can be computed using
\eqref{eq:DeltaR_lemma_P2} as
\begin{eqnarray}
&& \hspace{-15mm} E \left\{ \text{Tr} \left\{ \bV^{\dag} \Delta \bR
\bP^{\bot} d\bP \right\} \right\} = \frac{1}{N}\sum_{k = 1}^K
\frac{1}{2 j \left(
\ba_k^{(1)H} \bP^{\bot} \ba_k^{(1)} \right)} \nonumber \\
&& \hspace{30mm} \times \Bigg( \text{Tr} \Bigg\{ \bR \bP^{\bot}
\frac{\partial \bA}{\partial \omega_k} \left( \bA^H \bA \right)^{-1}
\bA^H \bV^{\dag} \bR \bV^{\dag}
\ba_k \ba_k^{(1)H} \bP^{\bot} \Bigg\} \nonumber \\
&& \hspace{30mm} - \text{Tr} \Bigg\{ \bR \bP^{\bot} \frac{\partial
\bA}{\partial \omega_k} \left( \bA^H \bA \right)^{-1} \bA^H
\bV^{\dag} \bR \bP^{\bot} \ba_k^{(1)} \ba_k^H \bV^{\dag} \Bigg\}
\Bigg). \label{eq:rho2_E_Tr_VdagDRPbotdP}
\end{eqnarray}
The second trace function in \eqref{eq:rho2_E_Tr_VdagDRPbotdP}
equals zero as $\bV^{\dag} \bR \bP^{\bot} = \boldsymbol{0}$. Then,
expression \eqref{eq:rho2_E_Tr_VdagDRPbotdP} can be rewritten as
\begin{eqnarray}
&& \hspace {-12mm} E \left\{ \text{Tr} \left\{ \bV^{\dag} \Delta \bR
\bP^{\bot} d\bP \right\} \right\} = \frac{\sigma_{\text{n}}^2}{N}
\sum_{k = 1}^K \frac{\ba_k^{(1)H} \bP^{\bot} \frac{\partial
\bA}{\partial \omega_k} \left( \bA^H \bA \right)^{-1} \bA^H
\bV^{\dag} \bR \bV^{\dag} \ba_k}{2 j \left( \ba_k^{(1)H} \bP^{\bot}
\ba_k^{(1)} \right)} \label{eq:rho2_E_Tr_VdagDRPbotdP_simplified}
\end{eqnarray}
where we used the equality $\bP^{\bot} \bR = \sigma_{\text{n}}^2
\bP^{\bot}$.

In a similar way, using expressions
\eqref{eq:dP_partialP_Deltaomega} and
\eqref{eq:partial_PA_partial_omegak}, the expected value of the
second trace function in \eqref{eq:rho2_final_simplified} is given
by
\begin{eqnarray}
&& \hspace{-15mm} E \left\{ \text{Tr} \left\{ d\bP \bP^{\bot} d\bP
\right\} \right\} = \nonumber \\
&& E \Bigg\{ \text{Tr} \Bigg\{ \sum_{k = 1}^K \sum_{i = 1}^K \bA
\left( \bA^H \bA \right)^{-1} \left( \frac{\partial \bA}{\partial
\omega_k} \right)^H \bP^{\bot} \frac{\partial \bA}{\partial
\omega_i} \left( \bA^H \bA \right)^{-1} \bA^H \Delta \omega_k \Delta
\omega_i\Bigg\} \Bigg\}. \label{eq:rho2_E_Tr_dPPbotdP_1}
\end{eqnarray}
Then, by substituting \eqref{eq:Delta_omega_k} in
\eqref{eq:rho2_E_Tr_dPPbotdP_1}, we have
\begin{eqnarray}
E \left\{ \text{Tr} \left\{ d\bP \bP^{\bot} d\bP \right\} \right\}
\hspace{-2mm} & = & \hspace{-2mm} E \Bigg\{ \text{Tr} \Bigg\{
\sum_{k = 1}^K \sum_{i = 1}^K \bA \left( \bA^H \bA \right)^{-1}
\left( \frac{\partial \bA}{\partial \omega_k} \right)^H \bP^{\bot}
\frac{\partial \bA}{\partial \omega_i} \left( \bA^H \bA \right)^{-1}
\bA^H \nonumber \\
&& \hspace{2mm} \times \frac{1}{2 j \left( \ba_k^{(1)H} \bP^{\bot}
\ba_k^{(1)} \right)} \times \frac{1}{2 j \left( \ba_i^{(1)H}
\bP^{\bot} \ba_i^{(1)} \right)} \nonumber \\
&& \hspace{2mm} \times \Big( \text{Tr} \left\{ \Delta \bR \bV^{\dag}
\ba_k \ba_k^{(1)H} \bP^{\bot} \right\} - \text{Tr} \left\{ \Delta
\bR \bP^{\bot} \ba_k^{(1)} \ba_k^H \bV^{\dag} \right\} \Big)
\nonumber \\
&& \hspace{2mm} \times \Big( \text{Tr} \left\{ \Delta \bR \bV^{\dag}
\ba_i \ba_i^{(1)H} \bP^{\bot} \right\} - \text{Tr} \left\{ \Delta
\bR \bP^{\bot} \ba_i^{(1)} \ba_i^H \bV^{\dag} \right\} \Big) \Bigg\}
\Bigg\} \nonumber \\
\end{eqnarray}
which is computed using \eqref{eq:DeltaR_lemma_P2} and the fact that
$\bP^{\bot} \bR \bV^{\dag} = \boldsymbol{0}$ as
\begin{eqnarray}
&& \hspace{-15mm} E \left\{ \text{Tr} \left\{ d\bP \bP^{\bot} d\bP
\right\} \right\} = \frac{\sigma_{\text{n}}^2}{2N} \sum_{k = 1}^K
\sum_{i = 1}^K \frac{\text{Tr} \left\{ \left( \frac{\partial
\bA}{\partial \omega_k} \right)^H \bP^{\bot} \frac{\partial
\bA}{\partial \omega_i} \left( \bA^H \bA \right)^{-1} \right\}}{
\left( \ba_k^{(1)H} \bP^{\bot} \ba_k^{(1)} \right) \left(
\ba_i^{(1)H} \bP^{\bot} \ba_i^{(1)} \right)} \nonumber \\
&& \hspace{52mm} \times Re \left\{ \ba_i^H \bV^{\dag} \bR \bV^{\dag}
\ba_k \ba_k^{(1)H} \bP^{\bot} \ba_i^{(1)} \right\}.
\label{eq:rho2_E_Tr_dPPbotdP}
\end{eqnarray}

Finally, the expected value of $\rho_2$ for a fixed value of
$\gamma$ is obtained using expressions
\eqref{eq:rho2_final_simplified},
\eqref{eq:rho2_E_Tr_VdagDRPbotdP_simplified}, and
\eqref{eq:rho2_E_Tr_dPPbotdP} as
\begin{eqnarray}
E \left\{ \rho_2 \right\} \hspace{-2mm} & = & \hspace{-2mm} \left( 1
- 2\gamma + \gamma^2 \right) E \left\{ \rho_1 \right\} \nonumber \\
&& \hspace{-2mm} + \frac{ 2 \left( \gamma - \gamma^2 \right)
\sigma_{\text{n}}^2 }{NK} Re \left\{ \sum_{k = 1}^K
\frac{\ba_k^{(1)H} \bP^{\bot} \frac{\partial \bA}{\partial \omega_k}
\left( \bA^H \bA \right)^{-1} \bA^H \bV^{\dag} \bR \bV^{\dag}
\ba_k}{2 j \left( \ba_k^{(1)H} \bP^{\bot} \ba_k^{(1)} \right)}
\right\} \nonumber \\
&& \hspace{-2mm} + \frac{\gamma^2 \sigma_{\text{n}}^2}{2NK} \sum_{k
= 1}^K \sum_{i = 1}^K \frac{\text{Tr} \left\{ \left( \frac{\partial
\bA}{\partial \omega_k} \right)^H \bP^{\bot} \frac{\partial
\bA}{\partial \omega_i} \left( \bA^H \bA \right)^{-1} \right\}}{
\left( \ba_k^{(1)H} \bP^{\bot} \ba_k^{(1)} \right) \left(
\ba_i^{(1)H} \bP^{\bot} \ba_i^{(1)} \right)} Re \left\{ \ba_i^H
\bV^{\dag} \bR \bV^{\dag} \ba_k \ba_k^{(1)H} \bP^{\bot} \ba_i^{(1)}
\right\}.\nonumber \\
\end{eqnarray}
It concludes the derivation.

\begin{comment}

\section*{Acknowledgment}

The authors would like to thank...
\end{comment}

\ifCLASSOPTIONcaptionsoff
  \newpage
\fi

\end{document}